\RequirePackage{ifpdf}
\ifpdf 
\documentclass[pdftex]{sigma}
\else
\documentclass{sigma}
\fi

\usepackage{xypic}
\numberwithin{equation}{section}

\newtheorem{Theorem}{Theorem}[section]
\newtheorem{Lemma}[Theorem]{Lemma}
\newtheorem{Proposition}[Theorem]{Proposition}
{\theoremstyle{definition}
\newtheorem{Remark}[Theorem]{Remark}
\newtheorem{Definition}[Theorem]{Definition}
}

\begin{document}

\allowdisplaybreaks

\renewcommand{\PaperNumber}{114}

\FirstPageHeading

\ShortArticleName{Projective Metrizability and Formal Integrability}

\ArticleName{Projective Metrizability and Formal Integrability}

\Author{Ioan BUCATARU~$^\dag$ and Zolt\'an MUZSNAY~$^\ddag$}

\AuthorNameForHeading{I.~Bucataru and Z.~Muzsnay}

\Address{$^\dag$~Faculty of Mathematics, Al.I.Cuza University, B-dul Carol 11, Iasi, 700506, Romania}
\EmailD{\href{mailto:bucataru@uaic.ro}{bucataru@uaic.ro}}
\URLaddressD{\url{http://www.math.uaic.ro/~bucataru/}}

\Address{$^\ddag$~Institute of Mathematics, University of Debrecen, H-4010 Debrecen, Pf. 12, Hungary}
\EmailD{\href{mailto:muzsnay@science.unideb.hu}{muzsnay@science.unideb.hu}}
\URLaddressD{\url{http://www.math.klte.hu/~muzsnay/}}

\ArticleDates{Received August 25, 2011, in f\/inal form December 08, 2011;  Published online December 12, 2011}

\Abstract{The projective metrizability problem can be
  formulated as follows: under what conditions the
  geodesics of a given spray coincide with the
  geodesics of some Finsler space, as oriented curves. In Theorem~\ref{thm:pm} we reformulate the projective
  metrizability problem for a~spray in terms of a f\/irst-order partial dif\/ferential
  operator $P_1$ and a set of
  algebraic conditions on semi-basic $1$-forms. We discuss the
  formal integrability of $P_1$ using two suf\/f\/icient conditions
  provided by Cartan--K\"ahler theorem. We prove in Theorem~\ref{thm:inv} that the symbol of~$P_1$ is involutive and hence one
  of the two conditions is always satisf\/ied. While discussing the
  second condition, in Theorem~\ref{thm:obstr1} we prove that there is only one obstruction to the formal integrability of~$P_1$, and this
  obstruction is due to the curvature tensor of the induced nonlinear
  connection. When the curvature obstruction is satisf\/ied, the projective
  metrizability problem reduces to the discussion of the algebraic
  conditions, which as we show are always satisf\/ied in the analytic case.  Based on these
  results, we recover all classes of sprays that are known to be projectively
  metrizable: f\/lat sprays, isotropic sprays, and arbitrary
  sprays on 1- and 2-dimensional manifolds. We provide examples of
  sprays that are projectively metrizable without being Finsler metrizable.}

\Keywords{sprays; projective metrizability; semi-basic forms; partial
  dif\/ferential operators; formal integrability}

\Classification{49N45; 58E30; 53C60; 58B20; 53C22}

\section{Introduction}

\looseness=1
The projective metrizability problem for a homogeneous system of
second-order ordinary dif\/fe\-ren\-tial equations, which can be identif\/ied
with a spray $S$, seeks for a
Finsler metric $F$ whose geodesics coincide with the geodesics of
the spray $S$, up to an orientation preserving
repara\-me\-te\-rization. For the case when $S$ is a f\/lat spray this problem
was f\/irst studied by Hamel~\cite{hamel03}  and it is known as the Finslerian version of
Hilbert's fourth problem~\cite{paiva05, crampin08, szilasi07}. In the
general case it was Rapcs\'ak~\cite{rapcsak62} who obtained, in
local coordinates, necessary and suf\/f\/icient conditions for the
projective metrizability problem of a spray. Global formulations for the
projective metrizability problem where obtained by Klein and Voutier~\cite{klein68}, and by Szilasi and Vattam\'any~\cite{szilasi02}. It has
been shown that this is an essential problem in various f\/ields of
biology and physics~\cite{aim93}.

\looseness=-1
The projective metrizability problem can be formulated as a particular case of the inverse
problem of the calculus of variations. We refer to \cite{anderso92, crampin81, krupkova07, morandi90, sarlet82} for various approaches of
the inverse problem of the calculus of variations. One of this
approaches seeks for the existence of a multiplier matrix that
satisf\/ies four Helmholtz conditions~\cite{krupkova07, sarlet82}. In~\cite{bucataru09},
these four Helmholtz conditions where reformulated in terms of a
semi-basic $1$-form. For the particular case of the projective
metrizability problem, it has been shown in \cite{bucataru09} that
only two of the four Helmholtz conditions are independent. In this
work we discuss the formal integrability of these two Helmholtz
conditions using two suf\/f\/icient conditions provided by Cartan--K\"ahler
theorem. The approach in this work follows the one developed in~\cite{muzsnay06} for studying the Finsler metrizability problem for a spray.

\looseness=1
In Section \ref{sec:prelim} we recall f\/irst some basic aspects of the Fr\"olicher--Nijenhuis
theory on a mani\-fold~$M$~\cite{frolicher56, KMS93}. Then, we use this
theory on $TM$ and apply it to the natural objects that live on the
tangent space: vertical distribution, Liouville vector f\/ield, and
semi-basic forms~\mbox{\cite{grifone72, grifone00, deleon89}}.

In Section~\ref{sec:pms} we use the geometric setting developed in the
previous section to reformulate the projective metrizability
problem. In Theorem~\ref{thm:pm} we obtain a set of necessary and
suf\/f\/icient conditions, for the projective metrizability problem of a
spray, which consists of a set of algebraic equations~\eqref{palgebric} and
a set of dif\/ferential equations~\eqref{cjh} on semi-basic forms. The
set of dif\/ferential equations determine a f\/irst-order partial dif\/ferential operator~$P_1$, called the
projective metrizability operator, which acts on semi-basic $1$-forms.

In Section~\ref{sec:fi} we discuss the formal integrability of the
projective metrizability operator $P_1,$ using two suf\/f\/icient
conditions provided by Cartan--K\"ahler theorem. Based on this theorem
and Theorems~\ref{thm:inv} and~\ref{thm:obstr1} we conclude that there is only one
obstruction to the formal integrability of~$P_1$. This obstruction is
expressed in terms of the curvature tensor of the nonlinear connection
induced by the spray. In this work we pay attention to various cases when
the obstruction condition is automatically satisf\/ied. Another possibility, which we
leave for further work, is to add this obstruction to the projective
metrizability operator and discuss the formal integrability of the new
operator. Using dif\/ferent techniques, an alternative expression of the obstruction
condition was obtain in~\cite[Theorem 4.9]{szilasi02}.

\looseness=1
In Section~\ref{sec:cpm} we discuss some classes of sprays for which the
curvature obstruction is automatically satisf\/ied: f\/lat sprays,
isotropic sprays, and arbitrary sprays on $2$-dimensional
manifolds. For each of these classes of sprays, the projective
metrizability problem reduces to the discussion of the algebraic
conditions~\eqref{palgebric}, which as we show are always
satisf\/ied in the analytic case. Although, for these classes, the projective
metrizability problem has been discussed before by some authors, our
approach in this work is dif\/ferent.
Using dif\/ferent methods, it was demonstrated
in~\cite{crampin07a} that f\/lat sprays are projectively metrizable. In~\cite{crampin07b} it has been shown that isotropic sprays are
projectively equivalent with f\/lat sprays, and hence are projectively metrizable.
On a~$2$-dimensional manifold it has
been shown by Matsumoto that every spray is projectively related to a
Finsler space~\cite{matsumoto95}, by extending the original discussion
of Darboux~\cite{darboux94} about second-order
dif\/ferential equations.

We use a spray on a $2$-dimensional, considered
by Anderson and Thompson in~\cite{anderso92}, and a~projectively f\/lat
spray of constant f\/lag curvature, considered by Yang in~\cite{yang11}, to provide examples of
projectively metrizable sprays that are not Finsler metrizable.

\section{Preliminaries} \label{sec:prelim}

In this section we present the dif\/ferential geometric tools
we need to formulate and study the projective metrizability problem.

\looseness=1
A systems of second-order ordinary dif\/ferential equations on a
manifold $M$ can be identif\/ied with a second-order vector f\/ield
that is called a semispray. To each semispray one can associate a
geometric apparatus very useful to obtain qualitative information
regarding: the variations of its geodesics, their stability, as
well as the inverse problem of the calculus of variations,~\cite{bcd10}.
A global formulation for the geometric apparatus one can associate
to a semispray is due to Grifone \cite{grifone72} and can be
obtained using the Fr\"olicher--Nijenhuis theo\-ry~\cite{frolicher56}.

\subsection[Fr\"olicher-Nijenhuis theory]{Fr\"olicher--Nijenhuis theory} \label{subsec:fn}

In this subsection we recall and extend some aspects of the
Fr\"olicher--Nijenhuis theory, which will be applied in the next
subsection to vector valued dif\/ferential forms on tangent bundles.
For the classic and modern formulations of Fr\"olicher--Nijenhuis theory we
refer to \cite{frolicher56, grifone72, grifone00, KMS93, deleon89}.

In this work $M$ is a real, $n$-dimensional, smooth manifold. We
denote by $C^{\infty}(M)$, the ring of smooth functions on $M$,
and by ${\mathfrak X}(M)$, the $C^{\infty}(M)$-module of vector
f\/ields on $M$. Consider $\Lambda(M)=\bigoplus_{k\in {\mathbb
N}}\Lambda^k(M)$ the graded algebra of dif\/ferential forms on $M$. We
denote by $S^k(M)$ the space of symmetric $(0,k)$ tensors on $M$.
We also write $\Psi(M)=\bigoplus_{k\in {\mathbb N}}\Psi^k(M)$ for the
graded algebra of vector-valued dif\/ferential forms on $M$.

For $L\in \Psi^l(M)$, a vector valued $l$-form, we consider
$\tau_L: \Lambda^1(M) \otimes \Lambda^k(M) \to \Lambda^{k+l}(M)$,
or $\tau_L: \Psi^1(M) \otimes \Psi^k(M) \to \Psi^{k+l}(M)$, the
\emph{alternating operator} def\/ined as
\begin{gather}
\label{taull}  (\tau_L B)(X_1,\dots ,X_{k+l}) = \frac{1}{k!l!} \sum_{\sigma \in S_{k+l}}
\varepsilon(\sigma) B(L(X_{\sigma(1)},\dots , X_{\sigma(l)}),
X_{\sigma(l+1)},\dots , X_{\sigma(l+k)}),
\end{gather}
where
$X_1,\dots , X_{k+l} \in {\mathfrak X}(M)$ and $S_{k+l}$ is the
permutation group of $\{1,\dots , k+l\}$.

The restriction of $\tau_L$ to $\Lambda^{k+1}(M)\subset
\Lambda^1(M) \otimes \Lambda^k(M)$, or the restriction to
$\Psi^{k+1}(M)$, is a~derivation of degree $(l-1)$ and it
coincides with the \emph{inner product}~$i_L$, see \cite{grifone00,
KMS93}.  Inner product~$i_L$
is trivial on $\Lambda^0(M)=C^{\infty}(M)$, or
$\Psi^0(M)=\mathfrak{X}(M)$, and hence it is a \emph{derivation of type
$i_*$} \cite{grifone00} or an \emph{algebraic derivation}
\cite{KMS93}. Since it satisf\/ies the Leibniz rule, $i_L$ is uniquely
determined by its action on $\Lambda^1(M)$, or $\Psi^1(M)$, when it is
given by $i_LB=B\circ L$. For the particular case when $l=1$ and $L=\operatorname{Id}$ we have that
$i_{\operatorname{Id}}B=kB$ for all $B\in \Lambda^k(M)$, or $B\in
\Psi^k(M)$.

For a linear connection $\nabla$ on $M$ consider
$d^{\nabla}:\Psi^k(M)\to \Psi^{k+1}(M)$ the covariant exterior
derivative, see \cite[\S~11.13]{KMS93}, given by
\begin{gather} \nonumber   d^{\nabla}B(X_1,\dots , X_{k+1}) =
\sum_{i=1}^{k+1}
(-1)^{i+1}\nabla_{X_i}B(X_1,\dots , \hat{X}_i, \dots , X_{l+1})  \\
\label{dnablaL}
\phantom{d^{\nabla}B(X_1,\dots , X_{k+1}) =}{}
+ \sum_{1\leq i<j\leq k+1} (-1)^{i+j}B([X_i,
X_j], X_1,\dots , \hat{X}_i,\dots , \hat{X}_j,\dots , X_{k+1}).\!\!\!
\end{gather} The exterior derivative $d: \Lambda^k(M)\to
\Lambda^{k+1}(M)$ satisf\/ies also formula \eqref{dnablaL} for $B\in
\Lambda^k(M)$. Therefore, we will use the notation $d^{\nabla}$ to refer to
both, the covariant exterior derivative, or the exterior
derivative. For the latter case $d=d^{\nabla}$ does not depend on
the linear connection $\nabla$.

For a vector valued $l$-form $L$, consider the commutator of the inner
product $i_L$ and the (covariant) exterior derivative $d$
($d^{\nabla}$). This dif\/ferential operator is
denoted by  $d_L: \Lambda^k(M) \to \Lambda^{k+l}(M)$ ($d^{\nabla}_L:
\Psi^k(M) \to \Psi^{k+l}(M)$), it is given by
\begin{gather} d^{\nabla}_L=i_L\circ d^{\nabla}+(-1)^{l}d^{\nabla}\circ i_L,
\label{dil}
\end{gather} it is a derivation of degree $l$, which is
called the \emph{$($covariant$)$ exterior derivative with respect to}~$L$. Derivation $d_L$ ($d^{\nabla}_L$) commutes with the exterior
derivative $d$ ($d^{\nabla}$) and hence it is a \emph{derivation of type}
$d_*$ \cite{grifone00} or a \emph{Lie derivation} \cite{KMS93}.
Since it satisf\/ies the Leibniz rule, $d_L$ ($d_{L}^{\nabla}$) is
uniquely determined by its action on $\Lambda^0(M)=C^{\infty}(M)$
($\Psi^0(M)={\mathfrak X}(M)$). For the particular case when $l=1$ and
$L=\operatorname{Id}$ we have that
$d^{\nabla}_{\operatorname{Id}}=d^{\nabla}$. Therefore, we obtain
$d^{\nabla}\operatorname{Id}=T$, where $T$ is the torsion of the linear
connection $\nabla$.

For two vector valued forms $L\in \Psi^l(M)$ and $K\in \Psi^k(M)$,
their Fr\"olicher--Nijenhuis bracket $[L,K]$ is a vector
valued $(k+l)$-form, def\/ined by
\begin{gather} \label{dlk}
d_{[L,K]}=d_L\circ d_K -(-1)^{kl}d_K\circ d_L.
\end{gather}

For a vector valued $l$-form $L$ and a linear connection $\nabla$
on $M$ we obtain a derivation of degree~$l$ given by ${\mathcal
D}_L=\tau_L \nabla$. Hence, ${\mathcal D}_L: \Lambda^k(M) \to
\Lambda^{k+l}(M)$ (${\mathcal D}_L: \Psi^k(M) \to
\Psi^{k+l}(M)$), acts on (vector-valued) $k$-forms as follows:
\begin{gather} \left({\mathcal
D_L}B\right)(X_1,\dots ,X_{k+l})  = \frac{1}{l!k!}
\sum_{\sigma\in S_{k+l}}
\varepsilon(\sigma)\big(\nabla_{L(X_{\sigma(1)},\dots ,X_{\sigma(l)})}B\big)\left(X_{\sigma(l+1)},\dots ,
X_{\sigma(k+l)}\right). \label{DLlk}
\end{gather}

For the particular case when $l=1$ and $L=\operatorname{Id}$, we
denote the corresponding derivation of degree $1$ by $\mathcal{D} =
\mathcal{D}_{\operatorname{Id}}$.
Since any derivation of degree $l$ can be uniquely
decomposed into a sum of a Lie derivation and an algebraic
derivation \cite[\S~8.3]{KMS93}, we obtain for ${\mathcal D}_L$
the following result.

\begin{Lemma} \label{lem:dec_DL} For a vector valued $l$-form $L$ and a linear connection $\nabla$ on $M$,
derivation of degree~$l$,~${\mathcal D}_L$, decomposes uniquely
into a sum of a Lie derivation and an algebraic derivation
as follows:
\begin{gather}
{\mathcal D}_L= d^{\nabla}_L -
i_{d^{\nabla}_L\operatorname{Id}}.\label{decompDL}
\end{gather}
\end{Lemma}
\begin{proof}
When acting on forms, formula \eqref{decompDL} reads ${\mathcal
D}_L= d_L - i_{d^{\nabla}_L\operatorname{Id}}$. The vector valued
$(l+1)$-form that def\/ines the inner product in formula
\eqref{decompDL} is, according to formula \eqref{dil}, given by
$d^{\nabla}_L\operatorname{Id}=i_LT +(-1)^ld^{\nabla}L.$

Since Lie derivations commute with the exterior
derivative $d^{\nabla}$ and satisfy the Leibnitz rule it follows that
they are uniquely determined by their action on
$\Lambda^0(M)=C^{\infty}(M)$ ($\Psi^0(M)={\mathfrak X}(M)$).
Using formulae \eqref{dil} and \eqref{DLlk} one can immediately
check that ${\mathcal D}_Lf=d_Lf$ for any scalar (vector
valued) $0$-form $f$.

Since algebraic derivations are trivial on
$\Lambda^0(M)=C^{\infty}(M)$ ($\Psi^0(M)={\mathfrak X}(M)$), and
satisfy the Leibnitz rule it follows that they are uniquely determined by
their action on $\Lambda^1(M)$ ($\Psi^1(M)$). To prove formula
\eqref{decompDL} we have to show that
\[
\left({\mathcal D}_L- d^{\nabla}_L\right)\omega = -\omega \circ
\left(d^{\nabla}_L\operatorname{Id}\right),
\]
 for any (vector valued) $1$-form
$\omega$. Since formally, we have the same formulae to def\/ine the
action on scalar, or vector valued forms, we will work with scalar
forms.

Let $\omega \in \Lambda^1(M)$ and $X_1, \dots , X_{l+1} \in
{\mathfrak X}(M)$. For $k=1$, from formula \eqref{DLlk} we obtain
the action of the derivation ${\mathcal D}_L$ on $1$-forms as
follows:
\begin{gather}
 \left({\mathcal D}_L\omega\right)(X_1,\dots , X_{l+1})=
\sum_{i=1}^{l+1}(-1)^{l+1-i}\left(\nabla_{L(X_1,\dots ,
\hat{X}_i,\dots ,X_{l+1})}\omega\right)(X_i) \nonumber  \\
\qquad{}
=  \sum_{i=1}^{l+1}(-1)^{l+1-i} \left\{L(X_1,\dots , \hat{X}_i,\dots ,X_{l+1})
(\omega(X_i)) -\omega\left( \nabla_{L(X_1,\dots ,
\hat{X}_i,\dots ,X_{l+1})} X_i\right)\right\}. \!\!\!\label{eq:dlo1}
\end{gather}

From formula \eqref{dil} we obtain that the action of the exterior
derivative $d_L$ on a $1$-form $\omega$ is given by
$d_L\omega=i_Ld\omega+ (-1)^ld(\omega\circ L)$.
Therefore, for $X_1, \dots , X_{l+1} \in {\mathfrak X}(M)$ we have{\samepage
\begin{gather}
\nonumber (d_L\omega)(X_1,\dots ,X_{l+1})   =   \sum_{i=1}^{l+1}(-1)^{l+1-i}
L(X_1,\dots , \hat{X}_i,\dots , X_{l+1})(\omega(X_i))  \\
\label{eq:dlo2}
\phantom{(d_L\omega)(X_1,\dots ,X_{l+1})   =}{}
+\sum_{i=1}^{l+1}(-1)^{l+1-i} \omega([X_i, L(X_1,\dots ,
\hat{X}_i,\dots , X_{l+1})])  \\
\nonumber   \phantom{(d_L\omega)(X_1,\dots ,X_{l+1})   =}{}
 +\sum_{1\leq i<j\leq l+1}
(-1)^{l+i+j}\omega(L([X_i,X_j], X_1,\dots , \hat{X}_i,\dots ,
\hat{X}_j,\dots , X_{l+1})).
\end{gather}
Now, we evaluate $d^{\nabla}_L\operatorname{Id}=i_LT + (-1)^l
d^{\nabla}L$ on $l+1$ vectors $X_1, \dots , X_{l+1} \in {\mathfrak
X}(M)$.}

For $k=1$, if we restrict the action of $\tau_L$ given by formula
\eqref{taull} to $\Psi^2(M)$ we obtain that the inner product $i_L: \Psi^2(M) \to
\Psi^{l+1}(M)$ is given by:
\begin{gather}
\left(i_LT\right)(X_1,\dots ,
X_{l+1})=\sum_{i=1}^{l+1}(-1)^{l-i}T(X_i, L(X_1,\dots , \hat{X}_i,
\dots , X_{l+1})). \label{eq:ilt}
\end{gather} Using formula \eqref{dnablaL}, the action of the exterior covariant derivative
$d^{\nabla}$ on the vector valued
$l$-form~$L$ is given by
\begin{gather}
\nonumber   d^{\nabla}L(X_1,\dots , X_{l+1}) = \sum_{i=1}^{l+1}
(-1)^{i+1}\nabla_{X_i}L(X_1,\dots , \hat{X}_i, \dots , X_{l+1})  \\
\phantom{d^{\nabla}L(X_1,\dots , X_{l+1}) =}{}
 + \sum_{1\leq i<j\leq l+1} (-1)^{i+j}L([X_i, X_j], X_1,\dots ,
\hat{X}_i,\dots , \hat{X}_j,\dots , X_{l+1}). \label{eq:dnl}
\end{gather}
Using formulae \eqref{eq:dlo1}, \eqref{eq:dlo2}, \eqref{eq:ilt}, and \eqref{eq:dnl} it follows that
\[
\left({\mathcal D}_L \omega-d_L\omega\right)(X_1,\dots ,X_{l+1}) =
-\left(\omega\circ (i_LT+(-1)^ld^{\nabla}L)\right)(X_1,\dots ,X_{l+1}),
\]
for all $X_1,\dots , X_{l+1}\in  {\mathfrak X}(M)$, which means that the decomposition \eqref{decompDL} is true.
\end{proof}

For the particular case when $l=1$ and $L=\operatorname{Id}$, we have that
the vector valued $2$-form
$d_{\operatorname{Id}}^{\nabla}\operatorname{Id}$ reduces to
torsion $T$ since $i_{\operatorname{Id}}T=2T$,
$d^{\nabla}\operatorname{Id}=T$. Therefore, decomposition
\eqref{decompDL} becomes
\[ 
{\mathcal D}=d- i_T.
\]
\begin{Remark}
Formula \eqref{decompDL} shows that the dif\/ference of the two
derivations $d_L-{\mathcal
D}_L=i_{d^{\nabla}_L\operatorname{Id}}$ is an algebraic derivation.
In other words, if $\omega \in \Lambda^k(M)$ vanishes at some
point $p\in M$, $\omega_p=0$, then $({\mathcal
D}_L\omega)_p=(d_L\omega)_p$. For the particular case when $l=1$, this result has been
shown in \cite[Proposition~2.5]{grifone00}.
\end{Remark}

\subsection[Vertical calculus on $TM$ and semi-basic forms]{Vertical calculus on $\boldsymbol{TM}$ and semi-basic forms} \label{subsec:vcalc}

Consider $(TM, \pi, M)$, the tangent bundle of the manifold $M$
and $(TM\setminus\{0\}), \pi, M)$ the slashed tangent bundle, which is
the tangent bundle with the zero
section removed. The tangent bundle carries some canonical structures, such
as the vertical distribution, the Liouville vector f\/ield, and the
vertical endomorphism. The dif\/ferential calculus associated to
these structures, using the Fr\"olicher--Nijenhuis theory developed in
the previous subsection, plays an important role in the geometry of a
system of second-order ordinary dif\/ferential equations,
\cite{ bcd10, bucataru09, grifone72, grifone00, klein68, deleon89, miron94}.

The \emph{vertical subbundle} is def\/ined as $VTM=\{\xi \in TTM,
(D\pi)(\xi)=0\}$. It induces a vertical distribution $V: u\in TM
\mapsto V_u=VTM\cap T_uTM$. This distribution is $n$-dimensional
and it is integrable, being tangent to the leaves of the \emph{natural
foliation} induced by submersion~$\pi$. If~$(x^i)$ are local coordinates on the base
manifold~$M$, we denote by $(x^i, y^i)$ the induced coordinates on~$TM$. It follows that~$y^i$ are coordinates along the leaves of the
natural foliation, while $x^i$ are transverse coordinates for the
foliation. We denote by ${\mathfrak X}^v(TM)$ the Lie subalgebra of vertical vector f\/ields on~$TM$. An
important vertical vector f\/ield on $TM$ is the \emph{Liouville vector
field}, which locally is given by
$\mathbb{C}=y^i{\partial}/{\partial y^i}.$

The \emph{tangent structure} (or vertical endomorphism) is the
$(1,1)$-type tensor f\/ield $J$ on~$TM$, which locally can be
written as follows:
\[
J=\frac{\partial}{\partial y^i} \otimes dx^i.
\]
Tensor $J$ satisf\/ies $J^2=0$ and $\operatorname{Ker} J=
\operatorname{Im} J=VTM$. Tangent structure $J$ is an integrable structure
since the Fr\"olicher--Nijenhuis bracket vanishes, $[J,J]=0$. As a consequence and using formula~\eqref{dlk} we have that
$d_J^2=0$.

For the natural foliation induced by submersion~$\pi$ and the
corresponding vertical distribution there are some important classes
of forms: basic and semi-basic forms. As we will see in the next
sections, semi-basic forms, vector valued semi-basic forms, and vector valued almost
semi-basic forms are important ingredients to formulate and
address the projective metrizability problem.
\begin{Definition} Consider $\omega\in \Lambda(TM)$ and $L\in \Psi(TM)$.
\begin{itemize}\itemsep=0pt
\item[i)] $\omega$ is called a \emph{basic} form if both $\omega$ and
  $d\omega$ vanish whenever one of the arguments of $\omega$
  (respectively $d\omega$) is a vertical vector f\/ield.
\item[ii)] $\omega$ is called a \emph{semi-basic} form if it vanishes whenever
one of its arguments is a vertical vector f\/ield.
\item[iii)] $L$ is called a vector valued \emph{semi-basic} form if it takes vertical values and vanishes
whenever one of its arguments is a vertical vector f\/ield.
\item[iv)] $L$ is called a vector valued \emph{almost semi-basic} form
  if it vanishes whenever one of its arguments is a vertical vector f\/ield
and for every vertical vector f\/ield $X\in {\mathfrak X}^v(TM)$ we
have that $\mathcal{L}_XL$ is a vector valued semi-basic form.
\end{itemize}
\end{Definition}
In local coordinates a basic $k$-form $\omega$ on $TM$ can
be written as
\begin{gather*}
 \omega=\frac{1}{k!}\omega_{i_1\dots i_k}(x)dx^{i_1}\wedge \cdots
\wedge dx^{i_k}.
\end{gather*}
For basic forms, the coordinates functions $\omega_{i_1\dots i_k}(x)$ are basic
functions, which means that they are constant along the leaves of the
natural foliation.

Locally, a semi-basic $k$-form $\omega$ on $TM$ can
be written as
\begin{gather}
 \omega=\frac{1}{k!}\omega_{i_1\dots i_k}(x,y)dx^{i_1}\wedge \cdots
\wedge dx^{i_k}.\label{eq:smf}
\end{gather}
We will denote by $\Lambda^k_v(TM)$ the set of
semi-basic $k$-forms on $TM$. A $1$-form $\omega$ on $TM$ is
semi-basic if and only if $i_J\omega=\omega \circ J=0$.

In local coordinates, a vector valued semi-basic $l$-form $L$ on
$TM$ can be written as
\begin{gather}
L=\frac{1}{l!}L^{\underline{j}}_{i_1\dots i_l}(x,y)\frac{\partial}{\partial y^j}\otimes dx^{i_1}\wedge \cdots
\wedge dx^{i_l}.\label{eq:vsmf}
\end{gather} In this work all contravariant or covariant
indices, of some tensorial coef\/f\/icients, that refer to vertical components will be underlined. We
will denote by $\Psi^l_v(TM)$ the set of vector valued semi-basic
$l$-forms on $TM$. A vector valued $1$-form $L$ on $TM$ is
semi-basic if and only if $J\circ L=0$ and $i_JL=L \circ J=0$. The
tangent structure $J$ is a vector valued semi-basic $1$-form.

Locally, a vector valued almost semi-basic $l$-form $L$ on $TM$ can be
expressed as
\begin{gather}
 L =\frac{1}{l!}L^{j}_{i_1\dots i_l}(x)\frac{\partial}{\partial
x^j}\otimes dx^{i_1}\wedge \cdots \wedge dx^{i_l}  +
\frac{1}{l!}L^{\underline{j}}_{i_1\dots i_l}(x,y)\frac{\partial}{\partial
y^j}\otimes dx^{i_1}\wedge \cdots \wedge
dx^{i_l}.\label{eq:vasmf}
\end{gather}
For a vector $X$ on $TM$ and a vector valued $l$-form $L$
on $TM$, the Fr\"olicher--Nijenhuis bracket $[X,L]$, def\/ined by formula
\eqref{dlk}, is a vector valued $l$-form on $TM$ given by:
\begin{gather*}
[X,L](X_1,\dots ,X_l)=[X,L(X_1,\dots ,X_l)] - \sum_{i=1}^lL([X,X_i],
X_1,\dots ,\hat{X}_i,\dots ,X_l),  
\end{gather*}
for $X_1, \dots , X_l$ vector f\/ields on $TM$. Using the above formula and the fact that the vertical distribution is integrable it
follows that vector valued semi-basic forms are also almost
semi-basic. This can be seen also from the local expressions
\eqref{eq:vsmf} and \eqref{eq:vasmf}.

Next two lemmas give a good motivation for considering the class of vector
valued almost semi-basic forms. We will also see in Section~\ref{sec:fi} that the partial
dif\/ferential operator we use to discuss the projective metrizability
problem is def\/ined in terms of some vector valued almost
semi-basic forms.

\begin{Lemma} \label{lem:dlomega} Let $L$ be a vector valued almost semi-basic $l$-form on~$TM$.
Then, the differential operator $d_L$ preserves semi-basic forms,
$d_L: \Lambda^k_v(TM) \to \Lambda^{k+l}_v(TM).$
\end{Lemma}

\begin{proof}
Consider a vector valued almost semi-basic $l$-form $L$, locally given
by formula \eqref{eq:vasmf}, and a semi-basic $k$-form $\omega$, locally given by formula
\eqref{eq:smf}. Since $d_L$ is a derivation of degree $l$, it follows
that the $(k+l)$-form $d_L\omega$ can be
expressed locally as follows
\begin{gather}
\nonumber d_L\omega  =  \frac{1}{k!}d_L(\omega_{i_1\cdots
  i_k})\wedge  dx^{i_1}\wedge \cdots \wedge dx^{i_k}  \\[-0.5ex]
 \phantom{d_L\omega  =}{}  +  \frac{1}{k!}
\sum^k_{j=1}(-1)^{(j-1)l} \omega_{i_1\cdots
  i_k} dx^{i_1}\wedge \cdots \wedge d_Ldx^{i_j} \wedge \cdots \wedge
dx^{i_k}. \label{eq:dlo}
\end{gather}
Using the assumption that $L$ is a vector valued almost semi-basic
form, we show that all terms in the right hand side of the above formula
are semi-basic forms. Since $L$ vanishes whenever one of its arguments
is a vertical vector f\/ield it follows that for a function $f\in C^{\infty}(TM)$, $d_Lf=i_Ldf=df\circ
L$ is a semi-basic $l$-form. Hence $d_L(\omega_{i_1\cdots i_k})$ are semi-basic $l$-forms.

We will prove now that $d_Ldx^{i_j}=(-1)^ldd_Lx^{i_j} $ are semi-basic
$(l+1)$-forms. Using the local expression \eqref{eq:vasmf} of $L$, we have
\begin{gather*} d_Lx^{i_j}=i_Ldx^{i_j} = dx^{i_j}\circ L =\frac{1}{l!}
L^{i_j}_{i_1\dots i_l}(x) dx^{i_1}\wedge \cdots \wedge
dx^{i_l},
\end{gather*}
which are basic $l$-forms. Therefore, $dd_Lx^{i_j}$ are basic and hence
semi-basic $(l+1)$-forms.
One can conclude now that all terms in the right hand side of formula~\eqref{eq:dlo} are
semi-basic forms and hence~$d_L\omega$ is a semi-basic~$(k+l)$-form.
\end{proof}

\begin{Lemma} \label{lem:nablaomega} Consider $\nabla$ a linear connection on $TM$ such
that $\nabla J=0$ and a vector valued almost semi-basic $l$-form
$L$ on $TM$. Then, the differential operator ${\mathcal D}_L$
preserves semi-basic forms, ${\mathcal D}_L: \Lambda^k_v(TM) \to
\Lambda^{k+l}_v(TM).$
\end{Lemma}

\begin{proof}
Using formula \eqref{decompDL} and Lemma \ref{lem:dlomega} we have
that the dif\/ferential operator ${\mathcal D}_L$ preserves semi-basic forms
if and only if the algebraic derivation of degree $l$,
$i_{d^{\nabla}_L\operatorname{Id}}$ preserves semi-basic forms. Hence,
we will complete the proof if we show that the vector valued
$(l+1)$-form $d^{\nabla}_L\operatorname{Id}=i_LT+(-1)^ld^{\nabla}_L$
takes vertical values whenever one of its arguments is a vertical
vector f\/ield. Here~$T$ is the torsion of the linear connection $\nabla$.

Using formulae \eqref{eq:ilt} and \eqref{eq:dnl} we have
\begin{gather}
\nonumber   (-1)^l(d^{\nabla}_L\operatorname{Id})(X_1,\dots ,X_{l+1}) =(-1)^l(i_LT+(-1)^ld^{\nabla}_L)(X_1,\dots ,X_{l+1})
\\[-0.5ex] \nonumber
\qquad{}  =  \sum_{i=1}^{l+1}(-1)^{i+1}
\left\{\nabla_{L(X_1,\dots ,\hat{X}_i,\dots ,X_{l+1})}X_i + [X_i,
  L(X_1,\dots ,\hat{X}_i,\dots ,X_{l+1})]\right\} \\[-0.5ex]
 \qquad\quad{}  +  \sum_{1\leq i<j\leq
  l+1} (-1)^{i+j}L([X_i,X_j], X_1,\dots , \hat{X}_i,\dots , \hat{X}_j,\dots ,
X_{l+1}). \label{eq:dnlid}
\end{gather}
We will show now that whenever one of the
$(l+1)$ arguments of $d^{\nabla}_L\operatorname{Id}$ is a vertical
vector f\/ield, then the right hand side of formula~\eqref{eq:dnlid} is
a vertical vector f\/ield. Using the fact that $d^{\nabla}_L\operatorname{Id}$ is a~vector valued $(l+1)$-form, we will discuss only the case when $X_1$ is a vertical
vector f\/ield. Since~$L$ vanishes whenever one of its arguments is a
vertical vector f\/ield, the nonzero vector f\/ield that remains from the
right hand side of formula~\eqref{eq:dnlid}, when $X_1$ is a vertical
vector f\/ield, is
\begin{gather}
\nabla_{L(X_2,\dots ,X_{l+1})}X_1+[X_1,L](X_2,\dots ,X_{l+1}). \label{eq:nl1}
\end{gather} The condition $\nabla J=0$ implies that the linear
connection $\nabla$ preserves the vertical distribution and since $X_1$ is a vertical
vector f\/ield it follows that $\nabla_{L(X_2,\dots,X_{l+1})}X_1 $ is a
vertical vector f\/ield as well. Since $L$ is a vector valued almost
semi-basic $l$-form and $X_1$ is a vertical
vector f\/ield it follows that $[X_1, L]$ is a vector valued semi-basic
form, therefore it takes values into the vertical distribution and
hence  $[X_1,L](X_2,\dots ,X_{l+1})$ is a vertical vector f\/ield. It
follows that the vector f\/ield in formula~\eqref{eq:nl1} is vertical
and hence we have completed the proof.
\end{proof}

\section{Projective metrizability problem of a spray}
\label{sec:pms}

A system of homogeneous second-order ordinary dif\/ferential equations
on a manifold~$M$, whose coef\/f\/icients do not depend explicitly on time, can be identif\/ied
with a special vector f\/ield on~$TM$ that is called a spray. In
this section we address the following question, known as the
projective metrizability problem: for a given spray $S$ f\/ind necessary
and suf\/f\/icient conditions for the existence of a Finsler function $F$
such that the geodesics of $S$ and the geodesics of $F$ coincide up to
an orientation preserving reparameterization. We obtain such necessary and
suf\/f\/icient conditions in Theorem \ref{thm:pm} and these conditions are expressed
in terms of semi-basic $1$-forms.

Particular aspects of the projective metrizability problem were studied more than a century ago by Hamel
\cite{hamel03}. The problem was formulated rigorously in 1960's by Rapcs\'ak
\cite{rapcsak62} and Klein and Voutier \cite{klein68}. Yet, the projective
metrizability problem is far from being solved, and in the last decade
it has been intensively studied \cite{paiva05, bucataru09,
  crampin07b, crampin08,  shen01, szilasi07,  szilasi02, yang11}.

\subsection{Spray, nonlinear connection, and curvature}
\label{subsec:sn}

In this subsection, we start with a spray $S$ and use the Fr\"olicher--Nijenhuis theory
to derive a dif\/ferential calculus on $TM\setminus\{0\}$ \cite{grifone72} and to
obtain information about the given system of SODE. For the remaining
part of the paper, all geometric objects will be considered def\/ined on
the slashed tangent bundle $TM\setminus\{0\}$ and not on the whole
$TM$. This is motivated by the fact that we will want to connect them
with geometric structures in Finsler geometry, where the Finsler
function is not dif\/ferentiable on the zero section.

\begin{Definition} \label{def:spray} A vector f\/ield $S\in{\mathfrak X}(TM\setminus\{0\})$
is called a \emph{spray} if
\begin{itemize}\itemsep=0pt
\item[i)] $JS=\mathbb{C}$,
\item[ii)] $[\mathbb{C}, S]=S$. \end{itemize}
\end{Definition}
First condition in Def\/inition~\ref{def:spray} expresses that a
spray $S$ can be locally given as
\[
 S=y^i\frac{\partial}{\partial x^i} - 2G^i(x,y) \frac{\partial}{\partial   y^i},
  \]
   for some functions $G^i$ def\/ined on domains of induced
coordinates on $TM\setminus\{0\}$.

Second condition in Def\/inition~\ref{def:spray} expresses that the
vector f\/ield $S$ is $2$-homogeneous. It is equivalent with the
fact that functions $G^i$ are $2$-homogeneous in the f\/ibre
coordinates. In this work we will consider positive homogeneity
only and hence $G^i(x,\lambda y)=\lambda^2 G^i(x,y)$
for all $\lambda>0$. By Euler's theorem this homogeneity condition is equivalent to
$\mathbb{C}(G^i)=2G^i$.

A curve $c:I \to M$ is called \emph{regular} if its tangent lift takes
values in the slashed tangent bundle, $c': I \to
TM\setminus\{0\}$. A regular curve is called a \emph{geodesic} of
spray $S$ if $S\circ c'=c''$. Locally, $c(t)=(x^i(t))$ is a
geodesic of spray $S$ if
\begin{gather*}
\frac{d^2x^i}{dt^2}+2G^i\left(x,\frac{dx}{dt}\right)=0.
\end{gather*}

\begin{Definition} A \emph{nonlinear connection} (or a horizontal
distribution, or Ehresmann connection) is def\/ined by an
$n$-dimensional distribution $H: u\in TM\setminus\{0\} \to
H_u\subset T_u(TM\setminus\{0\})$ that is supplementary to the
vertical distribution. \end{Definition}

Every spray induces a nonlinear connection through the
corresponding horizontal and vertical projectors, \cite{grifone72}
\begin{gather*}
h=\frac{1}{2}\left(\operatorname{Id}-\mathcal{L}_SJ\right), \qquad
v=\frac{1}{2}\left(\operatorname{Id}+\mathcal{L}_SJ\right).
\end{gather*} Locally, the two projectors $h$ and $v$ can
be expressed as follows
\[
h=\frac{\delta}{\delta x^i} \otimes dx^i, \qquad
v=\frac{\partial}{\partial y^i}\otimes \delta y^i,
\] where
\[
\frac{\delta}{\delta x^i}=\frac{\partial}{\partial
x^i}-N^j_i(x,y)\frac{\partial}{\partial y^j}, \qquad \delta y^i=dy^i+
N^i_j(x,y)dx^j, \qquad N^i_j(x,y)=\frac{\partial G^i}{\partial y^j}(x,y).
\]
Horizontal projector $h$ is a vector valued almost semi-basic
$1$-form.

For a spray $S$ consider the vector valued
semi-basic $1$-form
\begin{gather*}
\Phi=-v\circ \mathcal{L}_Sv = v\circ \mathcal{L}_Sh=v\circ
\mathcal{L}_S \circ h,
\end{gather*} which will be called the
\emph{Jacobi endomorphism}. It is also known as the Douglas tensor
\cite[Def\/inition~3.17]{grifone00} or as the Riemann curvature \cite[Def\/inition~8.1.2]{shen01}.
 Locally, the Jacobi endomorphism can be expressed as follows
\begin{gather*}
\Phi= R^i_j(x,y) \frac{\partial}{\partial y^i}
\otimes dx^j, \qquad R^i_j = 2\frac{\delta G^i}{\delta x^j} -
S(N^i_j) + N^i_kN^k_j. 
\end{gather*}

Another important geometric structure induced by a spray $S$ is
the \emph{curvature tensor} $R$. It is the vector valued semi-basic $2$-form
\begin{gather} R=\frac{1}{2}[h,h]=\frac{1}{2}R^i_{jk}\frac{\partial}{\partial y^i}\otimes dx^j\wedge dx^k. \label{curvature}
\end{gather}
Locally, the components of the curvature tensor, $R^i_{jk}$, are given by
\[
R^i_{jk}=\frac{\delta N^i_j}{\delta x^k} - \frac{\delta N^i_k}{\delta
x^j}.
\] Curvature tensor $R$ expresses the obstruction to the
integrability of the nonlinear connection.  Using formulae~\eqref{dlk}
and~\eqref{curvature} we have that $d_h^2=d_R$.

All the geometric objects induced by a spray $S$ inherit the
homogeneity condition. Therefore $[\mathbb{C}, h]=0$, which means
that the nonlinear connection is $1$-homogeneous. Also $[\mathbb{C},
R]=0$, $[\mathbb{C},\Phi]=\Phi$ and hence the the curvature tensor $R$
is $1$-homogeneous, while the  Jacobi endomorphism $\Phi$
is $2$-homogeneous.

Using the Jacobi identity, \cite[Proposition~2.7]{grifone00} , for the vector valued $0$-form $S$ and the
vector valued $1$-form $J$ we have $[J,[S,J]]-[J,[J,S]]-[S,[J,J]]=0$.
Therefore, we obtain $[J,h]=-2[J,[S,J]]=0$.

The two semi-basic vector vector valued $1$ and $2$-forms $\Phi$ and
$R$ are related as follows:
\begin{gather}
\Phi=i_SR, \qquad [J,\Phi]=3R. \label{eq:phir} \end{gather}
First formula in \eqref{eq:phir} is a consequence of the homogeneity,
while the second one is true in a~more general context. Locally, the
above two formulae can be expressed as follows:
\begin{gather*}
R^i_j=R^i_{kj}y^k, \qquad R^i_{jk}=\frac{1}{3}\left(\frac{\partial
    R^i_k}{\partial y^j}-\frac{\partial R^i_j}{\partial
    y^k}\right).
    \end{gather*}
An important class of sprays, which we will use in the last section to
provide examples of projectively metrizable sprays, is that of
isotropic sprays, \cite[Def\/inition~3.29]{grifone00}.
\begin{Definition} A spray $S$ is called \emph{isotropic} if its Jacobi
  endomorphism has the form
\begin{gather}
\Phi=\lambda J + \eta \otimes \mathbb{C}, \label{eq:phi_iso}
\end{gather}
where $\lambda\in C^{\infty}(TM\setminus\{0\})$ and $\eta$ is a
semi-basic $1$-form on $TM\setminus\{0\}$.
\end{Definition}
Due to f\/irst formula in \eqref{eq:phir} we have that $i_S\Phi=0$ and
hence $\lambda=-i_S\eta$. Also formulae \eqref{eq:phir} allows us to express
the isotropy condition \eqref{eq:phi_iso} for a spray in terms of the curvature tensor
$R$.
\begin{Proposition} \label{prop:iso}
A spray $S$ is isotropic if and only if its curvature tensor $R$ has the form
\begin{gather}
R=\alpha\wedge J + \beta \otimes \mathbb{C}, \label{eq:r_iso}
\end{gather}
where $\alpha$ is a semi-basic $1$-form and $\beta$ is a semi-basic
$2$-form on $TM\setminus\{0\}$. \end{Proposition}
\begin{proof}
We will prove that formulae \eqref{eq:phi_iso} and \eqref{eq:r_iso}
are equivalent.

Suppose that spray $S$ is isotropic. Therefore, the Jacobi endomorphism
$\Phi$ satisf\/ies formu\-la~\eqref{eq:phi_iso}. Using second formula
\eqref{eq:phir}, the formulae for the Fr\"olicher--Nijenhuis bracket of
two vector valued forms \cite[Appendix~A1]{grifone00}, and $[J, \mathbb{C}]=J$, we have
\begin{gather*}
3R=[J,\Phi]=[J, \lambda J + \eta\otimes \mathbb{C}] =
\left(d_J\lambda-\eta\right)\wedge J + d_J\eta \otimes \mathbb{C}. 
\end{gather*} Hence, the curvature tensor $R$ has the form
\eqref{eq:r_iso}.

We assume now that the curvature tensor $R$ has the form
\eqref{eq:r_iso}.  Using f\/irst formula \eqref{eq:phir} and the fact
that the inner product $i_S$ is a derivation of degree $-1$, we have
that the Jacobi endomorphism has the form
\begin{gather*}
\Phi=i_SR=i_S\alpha J + \left(i_S\beta-\alpha\right)\otimes
\mathbb{C}. 
\end{gather*} Hence, the spray $S$ is isotropic.
\end{proof}

We will use  Proposition \ref{prop:iso} and formula \eqref{eq:r_iso} in
Subsection \ref{subsec:pms} to show that isotropic sprays are
projectively metrizable sprays.

\subsection{Projectively related sprays} \label{subsec:prs}

Two sprays are projectively equivalent if their geodesics
coincide as oriented curves. Therefore, a spray is called projectively
metrizable if its geodesics coincide, as oriented curves, with the
geodesics of a Finsler space.

\looseness=-1
In \cite{bucataru09} it has been shown that the Helmholtz conditions for an
arbitrary semispray to be a~Lag\-rangian vector f\/ield can be
reformulated in terms of semi-basic $1$-forms. It has been shown also
that out of the four classic Helmholtz conditions only two of them are necessary and suf\/f\/icient
in the case of the projective metrizability problem for a spray. In this
subsection we obtain directly the two Helmholtz conditions, for
projective metrizability, in terms of semi-basic $1$-forms

\begin{Definition} \label{def:finsler}
By a \emph{Finsler function} we mean a continuous function $F: TM
\to \mathbb{R}$ satisfying the following conditions:
\begin{itemize}\itemsep=0pt
\item[i)] $F$ is smooth on $TM\setminus\{0\}$;
\item[ii)] $F$ is positive on $TM\setminus\{0\}$ and $F(x,0)=0$;
\item[iii)] $F$ is positively homogeneous of order $1$, which
means that $F(x,\lambda y)=\lambda F(x,y)$, for all $\lambda>0$
and $(x,y)\in TM$;
\item[iv)] the \emph{metric tensor} with components
\begin{gather*}
g_{ij}(x,y)=\frac{1}{2}\frac{\partial^2 F^2}{\partial y^i\partial
y^j} 
\end{gather*}
has rank $n$. \end{itemize}
\end{Definition}
According to Lovas \cite{lovas07}, conditions ii) and iv) of
Def\/inition~\ref{def:finsler}  imply that the metric tensor~$g_{ij}$ of
a Finsler function is positive def\/inite.

The regularity condition iv) of Def\/inition \ref{def:finsler} implies that the
Euler--Poincar\'e $2$-form of $F^2$, $\omega_{F^2}=dd_JF^2$, is non-degenerate
and hence it is a symplectic structure \cite{deleon89, morandi90}. Therefore, the equation
\begin{gather} i_Sdd_JF^2=-dF^2 \label{isddj}
\end{gather}
uniquely determine a vector f\/ield $S$ on $TM\setminus\{0\}$ that is called the
\emph{geodesic spray} of the Finsler function. Equation \eqref{isddj}
is equivalent to
\begin{gather} \mathcal{L}_Sd_JF^2=dF^2. \label{lsdj}
\end{gather}
Locally, the Euler--Poincar\'e $2$-form of $F^2$, $\omega_{F^2}=dd_JF^2$, can
be expressed as follows
\begin{gather*}
\omega_{F^2}=2g_{ij}\delta y^i \wedge dx^j. 
\end{gather*}
\begin{Definition} \label{defn:finsler_metr}
A spray $S$ is called \emph{Finsler metrizable} if there exists a
Finsler function $F$ that satisf\/ies one of the two equivalent
conditions \eqref{isddj} or \eqref{lsdj}.
\end{Definition}
One can reformulate condition iv) of Def\/inition \ref{def:finsler}
in terms of the Hessian of the Finsler function $F$ as follows.
Consider
\begin{gather*}
h_{ij}(x,y)=F\frac{\partial^2 F}{\partial y^i\partial y^j}
\end{gather*} the \emph{angular metric} of the Finsler function.
The metric tensor $g_{ij}$ and the angular tensor $h_{ij}$ are
related by
\begin{gather*}
g_{ij}=h_{ij} + \frac{\partial F}{\partial y^i} \frac{\partial
F}{\partial y^j}. 
\end{gather*} Metric tensor
$g_{ij}$ has rank $n$ if and only if angular tensor $h_{ij}$ has
rank $(n-1)$, see~\cite{matsumoto86}. Therefore, the regularity of
the Finsler function $F$ is equivalent with the fact that the
Euler--Poincar\'e $2$-form $\omega_{F}=dd_JF$ has rank $2n-2$.
\begin{Definition}\qquad
\begin{itemize} \itemsep=0pt \item[i)] Two sprays $S_1$ and $S_2$ are \emph{projectively equivalent}
if their geodesics coincide up to an orientation preserving
reparameterization. \item[ii)] A spray $S$ is \emph{projectively
metrizable} if it is projectively equivalent to the geodesic spray
of a Finsler function.
\end{itemize}
\end{Definition}
Two sprays $S_1$ and $S_2$ are projectively equivalent if and only
if there exists a $1$-homogeneous function $P\in
C^{\infty}(TM\setminus\{0\})$ such that $S_2=S_1-2P\mathbb{C}$,
\cite{aim93, shen01}.

Next theorem gives a characterization of projectively metrizable sprays in
terms of semi-basic $1$-forms on $TM\setminus\{0\}$.
\begin{Theorem} \label{thm:pm} A spray $S$ is projectively
metrizable if and only if there exists a semi-basic $1$-form
$\theta\in \Lambda^1_v(TM\setminus\{0\})$ such that
\begin{gather}
\operatorname{rank}\left(d\theta\right)=2n-2, \qquad i_S\theta>0,
\label{palgebric} \\ \mathcal{L}_{\mathbb{C}}\theta =0, \qquad
d_J\theta=0, \qquad d_h\theta=0. \label{cjh}\end{gather}
\end{Theorem}
\begin{proof}
We prove f\/irst that conditions \eqref{palgebric} and \eqref{cjh} are
necessary for the projective metri\-za\-bility problem of the spray $S$.
We assume that $S$ is projectively metrizable. Therefore, there
exists a~Finsler function $F$ with geodesic spray $S_F$ and a
$1$-homogeneous function $P$ on $TM\setminus\{0\}$ such that
$S=S_F-2P\mathbb{C}$. Consider $\theta=d_JF$, the Euler--Poincar\'e
$1$-form of the Finsler function~$F$. Due to the $1$-homogeneity
condition of $F$ it follows that $i_S\theta=\mathbb{C}(F)=F>0$. The non-degeneracy
of the Finsler function implies
$\operatorname{rank}\left(d\theta\right)=2n-2$. Since $\theta$ is
$0$-homogeneous it follows that
$\mathcal{L}_{\mathbb{C}}\theta=0$. Condition $d_J\theta=0$ is
also satisf\/ied since $d_J\theta=d_J^2F=0$.

\looseness=-1
It remains to show that $d_h\theta=0$. The geodesic spray $S_F$ is
uniquely determined by condi\-tion~\eqref{isddj}, from which it
follows that $S_F(F^2)=0$ and hence $S_F(F)=0$. Since $S_F$ also
satis\-f\/ies condition \eqref{lsdj} it follows that
$\mathcal{L}_{S_F}\left(F\theta\right)=FdF$, which implies
$\mathcal{L}_{S_F}\theta=dF$. Using $S=S_F-2P\mathbb{C}$ we obtain
that $\mathcal{L}_{S}\theta -2 \mathcal{L}_{P\mathbb{C}}\theta=dF$.
Using again the $0$-homogeneity of the semi-basic $1$-form $\theta$ it
follows $\mathcal{L}_{P\mathbb{C}}\theta= P
\mathcal{L}_{\mathbb{C}}\theta=0$ and hence
$\mathcal{L}_{S}\theta=dF$. We apply now $d_J$ to both sides of
this last relation and use the commutation rules
$\mathcal{L}_Sd_J-d_J\mathcal{L}_S=d_{[S,J]}=-d_{h}+d_{v}$ and
$dd_J+d_Jd=0.$ Therefore,
\[
-d_h\theta-d_v\theta=-dd_JF=d_JdF=d_J\mathcal{L}_{S}\theta=
\mathcal{L}_Sd_J\theta + d_{h}\theta- d_{v}\theta,
\] from where
it follows that $d_h\theta=0$.

We prove now that the conditions \eqref{palgebric} and \eqref{cjh}  are
suf\/f\/icient for the projective metrizability problem of the spray $S$. Consider $\theta\in \Lambda^1(TM\setminus\{0\})$
a semi-basic $1$-form that sa\-tis\-f\/ies conditions~\eqref{palgebric}
and~\eqref{cjh}. Def\/ine the function $F=i_S\theta$. Using the commutation rule
$i_Sd_J+d_Ji_S=\mathcal{L}_{\mathbb{C}}-i_{[S,J]}$ as well as
conditions $d_J\theta=0$ and $\mathcal{L}_{\mathbb{C}}\theta=0$ it
follows that $d_JF=d_Ji_S\theta=i_h\theta=\theta$. Hence $\theta$
is the Euler--Poincar\'e $1$-form of $F$. Now conditions~\eqref{palgebric} assure that $F$ is a Finsler function. Consider
the function $P\in C^{\infty}(TM\setminus\{0\})$ given by
$2P=S(F)/F$, which is $1$-homogeneous. We will show now that the
spray $\tilde{S}=S-2P\mathbb{C}$ satisf\/ies equation~\eqref{lsdj}
and hence it is the geodesic spray of the Finsler function~$F$.

Using the commutation rule $i_Sd_h+d_hi_S=\mathcal{L}_S-
i_{[S,h]}$ and the fact that $d_h\theta=0$ it follows
$0=i_Sd_h\theta=-d_hi_S\theta+ \mathcal{L}_S\theta-
i_{[S,h]}\theta$. Using the fact that $i_{[S,h]}\theta=dF\circ
J\circ \mathcal{L}_Sh=dF\circ v=d_vF$ it follows that
$\mathcal{L}_S\theta=d_hF+d_vF=dF$. We show now that $\tilde{S}$
satisf\/ies the same equation. Indeed
$\mathcal{L}_{\tilde{S}}\theta=
\mathcal{L}_{S-2P\mathbb{C}}\theta=dF$ since
$\mathcal{L}_{P\mathbb{C}}\theta=0$. From the def\/ining formula of
function $P$ is follows that
$\tilde{S}(F)=S(F)-2P\mathbb{C}(F)=S(F)-2PF=0$. Therefore
$\mathcal{L}_{\tilde{S}}d_JF^2=2F\mathcal{L}_{\tilde{S}}d_JF=2FdF=dF^2$
and hence $\tilde{S}$ is the geodesic spray of the Finsler
function~$F$.
\end{proof}

The second part of the proof of Theorem \ref{thm:pm} shows that if
there exists a semi-basic $1$-form $\theta$ on $TM\setminus\{0\}$ that
satisf\/ies the conditions \eqref{palgebric} and \eqref{cjh} then the
given spray $S$ is projectively related to the spray
\[
 S_F=S-\frac{\mathcal{L}_{S}(i_S\theta)}{i_S\theta}\mathbb{C},
 \]
  which is
the geodesic spray of the Finsler function $F=i_S\theta$. In this
case, the semi-basic $1$-form $\theta=\theta_i dx^i$ is the
Euler--Poincar\'e $1$-form of the Finsler function $F$,
$\theta=d_JF$. Therefore,
\begin{gather} \theta_i=\frac{\partial F}{\partial y^i}, \qquad
  h_{ij}=F\frac{\partial \theta_i}{\partial
    y^j}, \qquad Fd\theta=h_{ij} \delta y^i \wedge
  dx^j. \label{fdt} \end{gather}
Formulae \eqref{fdt} show the relation between a semi-basic $1$-form
$\theta$, a solution of the projective metrizability problem using
Theorem~\ref{thm:pm}, and the classic approach of the problem using the
multiplier matrix $h_{ij}$.

\section{Formal integrability for the projective metrizability problem}
\label{sec:fi}

Theorem \ref{thm:pm} provides necessary and suf\/f\/icient conditions for the projective
metrizability problem. These conditions consist of a set of algebraic equations~\eqref{palgebric},
and a set of dif\/ferential equations~\eqref{cjh}. In this section,
we study the set of dif\/ferential equations~\eqref{cjh} using Spencer's
technique of formal integrability \cite{BCG91, grifone00}  and following some of the techniques
used for studying the Finsler metrizability problem, which were developed in~\cite{muzsnay06}.

\subsection{Formal integrability} \label{subsec:fi}
In this subsection, we recall f\/irst the basic notions of formal integrability
\cite{BCG91, grifone00} and then we apply it to the system \eqref{cjh}.

Consider $E$ a vector bundle over the base manifold $M$. For a
section $s$ of $E$ and $k\geq 1$ we denote by $j^k_xs$ the $k$th
order jet of $s$ at the base point~$x$ in~$M$. The bundle of $k$th
order jets of sections of $E$ is denoted by~$J^kE$. For two vector
bundles~$E$ and~$F$ over the same base manifold~$M$, a linear
partial dif\/ferential operator of order~$k$,
\[
P: \  \operatorname{Sec}(E) \to \operatorname{Sec}(F),
\]  can be identif\/ied with a morphism of vector
bundles over $M$, $p^{0}(P): J^{k}E \to
F$. We will also consider the $l$th order jet prolongation of the dif\/ferential
operator $P$, which will be identif\/ied with the morphisms of vector bundles
over $M$, $p^{l}(P): J^{k+l}E \to J^lF$, def\/ined by
\[
p^{l}(P)\big(j^{k+l}_xs\big)=j^l_x(Ps).
\] We will denote by
$R^{k+l}_x(P)=\operatorname{Ker} p_x^{l}(P)\subset J^{k+l}_xE$ the
space of $(k+l)$th order formal solutions of $P$ at $x$ in $M$.
\begin{Definition} \label{def:fi} The dif\/ferential operator $P$ is called \emph{formally integrable} at
$x$ in $M$ if $R^{k+l}(P)$ is a vector bundle over $M$, for all $l\geq 0$,
and the map $\bar{\pi}^{k+l-1}_x: R^{k+l}_x(P) \to R^{k+l-1}_x(P)$ is onto for all
$l\geq 1$. \end{Definition}
In the analytic case, formal integrability implies
existence of analytic solutions for arbitrary initial data, see
\cite[p.~397]{BCG91}.

Denote by $\sigma^k(P): S^k(M)\otimes E \to F$ the symbol of $P$,
which is def\/ined by the highest order terms of the dif\/ferential operator
$P$, and by $\sigma^{k+l}(P): S^{k+l}(M)\otimes E \to S^l(M)
\otimes F$ the symbol of the $l$th order prolongation of $P$. For
each $x$ in $M$, we write
\begin{gather*}
g^k_x(P) = \operatorname{Ker}\sigma^k_x(P),
\\ g^k_x(P)_{e_1\dots e_j}  =  \{A \in g^k_x(P)| i_{e_1}A=\cdots =
i_{e_j}A=0\},\qquad j\in \{1,\dots ,n\},
\end{gather*}
where  $\{e_1,\dots ,e_n\}$ is a basis of $T_xM$. Such a basis is
called \emph{quasi-regular} if it satisf\/ies
\begin{gather} \dim g^{k+1}_x(P)= \dim g^k_x(P) + \sum_{j=1}^n \dim
g^k_x(P)_{e_1\dots e_j}. \label{quasir}\end{gather} The symbol
$\sigma^k(P)$ is called \emph{involutive} at
$x$ in $M$ if there exists a quasi-regular basis of $T_xM$.

In this work we will address the projective metrizability problem by
discussing f\/irst the formal integrability of the system \eqref{cjh}. For
this we will use the two suf\/f\/icient conditions provided by
Cartan--K\"ahler theorem.

\medskip

\noindent {\bf Theorem} [Cartan--K\"ahler]. {\it Let $P$ be a linear partial
differential operator of order $k$. Suppose $g^{k+1}(P)$ is a
vector bundle over $R^k(P)$. If the map $\overline{\pi}^k:
R^{k+1}(P) \to R^k(P)$ is onto and the symbol $\sigma^k(P)$ is involutive, then
$P$ is formally integrable.}

\medskip

In order to study the formal integrability of the system
\eqref{cjh} we consider the f\/irst-order partial dif\/ferential operator $P_1:
\Lambda^1_v(TM\setminus\{0\}) \to \Lambda^1_v(TM\setminus\{0\})
\oplus \Lambda^2_v(TM\setminus\{0\}) \oplus
\Lambda^2_v(TM\setminus\{0\})$, which we call \emph{the projective
  metrizability operator}
\begin{gather}
\label{p1} P_1  = \left({\mathcal L}_{\mathbb C}, d_J, d_h
\right).
\end{gather}
Since $\mathbb{C} $ and $J$ are vector valued, semi-basic $0$ and
respectively $1$-forms and
$h$ is a vector valued almost semi-basic $1$-form,  according to Lemma \ref{lem:dlomega}, all dif\/ferential
operators $\mathcal{L}_{\mathbb{C}}$, $d_J$, $d_h$ preserve
semi-basic forms. Therefore, the dif\/ferential operator $P_1$ is well
def\/ined.

\subsection{Involutivity of the projective metrizability operator}
\label{subsec:inv}

In this subsection we prove that the projective metrizability operator
\eqref{p1} satisf\/ies one of the two suf\/f\/icient conditions for formal
integrability, provided by Cartan--K\"ahler theorem, namely we will
prove that the symbol $\sigma^1(P_1)$ is involutive.

Since all the bundles we will refer to in this subsection are vector
bundles over $TM\setminus\{0\}$, we will omit mentioning it
explicitly. For example, we will denote by $T^*_v$ the vector
bundle of semi-basic $1$-forms $T^*_v(TM\setminus\{0\})$, which is
a subbundle of $T^*(TM\setminus\{0\})$. We will denote by
$\Lambda^kT^*_v$ the vector bundle of semi-basic $k$-forms on
$TM\setminus\{0\}$, and by
$\Lambda^k_v=Sec\left(\Lambda^kT^*_v\right)$ the
$C^{\infty}(TM\setminus\{0\})$-module of sections $\Lambda^k_v(TM\setminus\{0\})$. By $S^kT^*$ we denote the vector
bundle of symmetric tensors of $(0,k)$-type on $TM\setminus\{0\}$.

The partial dif\/ferential operator $P_1$ induces a morphism of vector
bundles
\[
 p^0(P_1): \ J^1T^*_v \to F_1:=T^*_v \oplus \Lambda^2 T^*_v \oplus
\Lambda^2 T^*_v.
\] Together with this morphism we will consider the $l$th order jet prolongations $p^l(P_1):
J^{l+1}T^*_v \to J^{l}F_1$, for $l\geq 1$.

Locally, for a semi-basic $1$-form $\theta=\theta_idx^i\in
\Lambda^1_v$, we have
\begin{gather*}
\mathcal{L}_{\mathbb C}\theta=\frac{\partial \theta_i}{\partial
  y^j}y^jdx^i, \qquad d_J\theta = \frac{1}{2}\left(\frac{\partial
\theta_i}{\partial y^j}-\frac{\partial \theta_j}{\partial
y^i}\right)dx^j\wedge dx^i, \\ d_h\theta = \frac{1}{2}\left(\frac{\delta \theta_i}{\delta x^j}-
\frac{\delta \theta_j}{\delta x^i} \right) dx^j\wedge dx^i.
\end{gather*}
Therefore, the vector bundle morphism $p^0(P_1)$ can be expressed as
follows
\begin{gather*}
p^0(P_1)(j^1\theta)= \left(\frac{\partial
\theta_i}{\partial y^j}y^j dx^i, \frac{1}{2}\left(\frac{\partial
\theta_i}{\partial y^j}-\frac{\partial \theta_j}{\partial
y^i}\right) dx^j\wedge dx^i, \frac{1}{2}\left(\frac{\delta \theta_i}{\delta x^j}-
\frac{\delta \theta_j}{\delta x^i}\right) dx^j\wedge dx^i \right).
\end{gather*} The symbol
of $P_1$ is the vector bundle morphism $\sigma^1(P_1): T^*\otimes
T^*_v \to F_1$, def\/ined by the highest order terms of $p^0(P_1)$.
Since all terms that def\/ine $p^0(P_1)$ are f\/irst-order terms,
it follows that
\begin{gather} \sigma^1(P_1)A=\left(
\sigma^1\left(\mathcal{L}_{\mathbb C}\right)A=i_{\mathbb{C}}A,
\sigma^1\left(d_J\right)A=\tau_JA,
\sigma^1\left(d_h\right)A=\tau_hA\right). \label{sigma1p1}
\end{gather}
In view of formula \eqref{taull}, the three components of the vector
bundle morphism $\sigma^1(P_1)$ are given~by:
\begin{gather*}
\left( \sigma^1\left(\mathcal{L}_{\mathbb C}\right)A\right) (X)
 = \left( i_{\mathbb{C}}A\right)(X)=A(\mathbb{C}, X); \\
\left(\sigma^1\left(d_J\right)A\right)(X,Y) = \left(\tau_JA\right)(X,Y)
= A(JX,Y)-A(JY, X); \\
\left(\sigma^1\left(d_h\right)A\right)(X,Y) = \left(\tau_hA\right)(X,Y)
= A(hX,Y)-A(hY, X),
\end{gather*} for $X$, $Y$ vector f\/ields on $TM\setminus\{0\}$.
Note that for $A\in T^*\otimes T^*_v$, $i_{\mathbb{C}}A$, $\tau_JA$,
$\tau_hA$ are semi-basic forms and hence the symbol $\sigma^1(P_1)$ is well
def\/ined.

The f\/irst-order prolongation of the symbol of $P_1$ is the vector
bundle morphism $\sigma^2(P_1): S^2T^*\otimes T^*_v \to T^*\otimes
F_1$ that satisf\/ies
$i_X\left(\sigma^2(P_1)B\right)=\sigma^1(P_1)(i_XB)$ for all $B\in
S^2T^*\otimes T^*_v$ and all $X\in {\mathfrak
X}(TM\setminus\{0\})$. Therefore, we obtain
\begin{gather*}
\sigma^2(P_1)B=\left( \sigma^2\left(\mathcal{L}_{\mathbb
C}\right)B, \sigma^2\left(d_J\right)B,
\sigma^2\left(d_h\right)B\right),
\end{gather*} where  for $X, Y, Z$ vector f\/ields on
$TM\setminus\{0\}$ we have:
\begin{gather} \nonumber \left(\sigma^2\left(\mathcal{L}_{\mathbb
C}\right)B\right)(X,Y)   =   B(X, \mathbb{C}, Y), \\
\label{sigma2p1} \left(\sigma^2\left(d_J\right)B\right)(X,Y,Z)   =
B(X, JY,Z)-B(X,JZ,Y), \\
\nonumber \left(\sigma^2\left(d_h\right)B\right)(X,Y,Z)   =
B(X,hY,Z)-B(X,hZ, Y).
\end{gather}
\begin{Theorem} \label{thm:inv}
The symbol $\sigma^1(P_1)$, of the projective metrizability operator $P_1=\left({\mathcal L}_{\mathbb C},
d_J, d_h \right)$, is involutive.
\end{Theorem}
\begin{proof}
The symbol $\sigma^1(P_1)$ is involutive if there exists
a quasi-regular basis of $T_u(TM\setminus\{0\})$. It means that we will
have to seek for a basis of $T_u(TM\setminus\{0\})$ that satisf\/ies
the equality \eqref{quasir} for $k=1$, at some point $u\in TM\setminus\{0\}$.

We start by computing the f\/irst term in the right hand side of formula
\eqref{quasir}, which is $\dim g^1_u(P_1)$, for some $u\in
TM\setminus\{0\}$. Recall that $g^1(P_1)=\operatorname{Ker} \left(\sigma^1(P_1)\right)
\subset T^*\otimes T^*_v$.  We have to compute the dimension of
the f\/ibers of $g^1(P_1)$, which is a vector subbundle of $ T^*\otimes T^*_v$. An element $A\in
g^1(P_1)$ can be expressed, with respect to the adapted dual basis
$\{dx^i, \delta y^i\}$, as follows
\begin{gather*}
A=A_{ij}dx^i\otimes dx^j + A_{\underline{i}j}\delta y^i \otimes
dx^j. 
\end{gather*} Using formula \eqref{sigma1p1},
the symbol $\sigma^1(P_1)$ can be expressed as follows:
\[
\sigma^1(P_1)A=\left( A_{\underline{i}j}y^idx^j, \frac{1}{2}(
A_{\underline{i}j}- A_{\underline{j}i})dx^i\wedge dx^j,
\frac{1}{2}(A_{ij}- A_{ji})dx^i\wedge dx^j\right).
\]
The condition $\tau_hA=0$ is equivalent with $A_{ij}=A_{ji}$ and due
to this condition $A_{ij}$ contribute with $n(n+1)/2$ to the $\dim g^1_u(P_1)$. The conditions
$\tau_JA=0$ and $i_{\mathbb{C}}A=0$ are equivalent to
$A_{\underline{i}j}=A_{\underline{j}i}$, and respectively
$A_{\underline{i}j}y^i=0$. Hence, due to these two conditions,
$A_{\underline{i}j}$ contribute
with $n(n-1)/2$ to the the $\dim g^1_u(P_1)$. It follows that $\dim
g^{1}_u(P_1)=n(n-1)/2 + n(n+1)/2 = n^2.$

We continue the proof by computing the left hand side of formula
\eqref{quasir}, which is $\dim g^2_u(P_1)$.
Therefore, we will consider the kernel of the f\/irst-order prolongation of the symbol,
$g^2(P_1)=\operatorname{Ker}\left(\sigma^2(P_1)\right)\subset
S^2T^*\otimes T^*_v$. An element $B\in g^2(P_1)$ can be expressed,
with respect to the adapted dual basis $\{dx^i, \delta y^i\}$, as
follows
\begin{gather*}
B=B_{ijk}dx^i\otimes dx^j\otimes dx^k + B_{\underline{i}jk}\delta
y^i \otimes dx^j \otimes dx^k \\
\phantom{B=}{} + B_{i\underline{j}k}dx^i \otimes
\delta y^j \otimes dx^k + B_{\underline{ij}k}\delta y^i \otimes
\delta y^j \otimes dx^k,
\end{gather*}
with the symmetry conditions $B_{ijk}=B_{jik}$,
$B_{\underline{i}jk}=B_{i\underline{j}k}$ and
$B_{\underline{ij}k}=B_{\underline{ji}k}$ satisf\/ied. Using formula~\eqref{sigma2p1}, the symbol $\sigma^2(P_1)$ can be expressed as
follows
\begin{gather*}
\sigma^2(P_1)B  = \bigg(B_{i\underline{j}k}y^j dx^i \otimes dx^k +
B_{\underline{ij}k}y^j \delta y^i \otimes dx^k,  \\
\hphantom{\sigma^2(P_1)B  = \bigg(}{}  \frac{1}{2}(B_{i\underline{j}k} -
B_{i\underline{k}j})dx^i\otimes dx^j \wedge dx^k +
\frac{1}{2}(B_{\underline{ij}k} - B_{\underline{ik}j})\delta
y^i\otimes dx^j \wedge dx^k, \\
\hphantom{\sigma^2(P_1)B  = \bigg(}{}  \frac{1}{2}(B_{ijk} - B_{ikj})dx^i\otimes dx^j \wedge
dx^k + \frac{1}{2}(B_{\underline{i}jk} -
B_{\underline{i}kj})\delta y^i\otimes dx^j \wedge dx^k\bigg).
\end{gather*}
The totally symmetric components $B_{ijk}$ contribute with
$n(n+1)(n+2)/6$ to the $\dim g^2_u(P_1)$. The other two are also
totally symmetric components on the $(n-1)$-dimensional space
given by restrictions $B_{i\underline{j}k}y^j=0$ and respectively
$B_{\underline{ij}k}y^j=0$. Therefore, each of them contributes
with $(n-1)n(n+1)/6$ to the $\dim g^2_u(P_1)$. Consequently, $\dim
g^2_u(P_1)=n(n+1)(n+2)/6+2(n-1)n(n+1)/6=n^2(n+1)/2$.

Finally, for some $u\in TM\setminus\{0\} $, we seek for a
basis of $T_u(TM\setminus\{0\})$ for
which formula \eqref{quasir} holds true.  Consider $\{h_i, i\in \{1,\dots ,n\}\}$ a basis for the horizontal
distribution and $\{v_i, i\in \{1,\dots ,n\}\}$, with $v_n=\mathbb{C}$,
a basis for the vertical distribution such that $Jh_i=v_i$, for all
$i\in \{1,\dots ,n\}$.

For $A\in g^1(P_1)$, and the basis $\mathcal{B}=\{h_i, v_i, i\in
\{1,\dots ,n\}\}$, let us denote
\begin{gather*}
a_{ij}=A(h_i, h_j) \qquad \mathrm{and} \qquad  b_{ij}=A(h_i,
v_j).\end{gather*}
It follows that \begin{enumerate}\itemsep=0pt
\item $a_{ij} = a_{ji}$, $i,j = 1,\dots ,n$, because $A \in
  \mathrm{Ker} \, \sigma^1(d_h)$,
\item $b_{ij} = b_{ji}$,   $i,j = 1,\dots ,n$, because $A \in
  \mathrm{Ker} \, \sigma^1(d_J)$,
\item $b_{ni} =(b_{in} )= 0$,   $i = 1,\dots ,n$, because $A \in
  \mathrm{Ker} \, \sigma^1(\mathcal L_{\mathbb{C}})$.
\end{enumerate}
Note that $\dim  g^1(P_1)=n^2$ is determined by the $n(n+1)/2$ independent components
  of $a_{ij}$ and $n(n-1)/2$ independent components of $b_{ij}$.
We will prove now that $\widetilde{\mathcal{B}}=\{e_i, v_i, i\in \{1,\dots ,n\}\}$,
where $e_1=h_1, e_2=h_2+v_1,\dots , e_{n-1}=h_{n-1}+v_{n-2},
e_n=S+v_{n-1}$, is a quasi-regular basis.

For the basis $\widetilde{\mathcal{B}}$ we denote
\begin{gather*}
  \tilde{a}_{ij}= A(e_i,e_j)  \qquad \mathrm{and} \qquad  \tilde{b}_{ij} = A(v_i,e_j).
\end{gather*}
Because $A$ is semi-basic in the second variable we have
\begin{gather*}
    A(e_i,e_j) = A(h_i+v_{i-1},h_j+v_{j-1}) = A(h_i+v_{i-1},h_j),
  \\
    A(v_i,e_j) = A(v_i,h_j+v_{j-1}) = A(v_i,h_j),
\end{gather*}
which means that
\begin{gather*}
  \tilde{a}_{ij}   = a_{ij} + b_{i-1,j}, \qquad   
   \tilde{b}_{ij}  = b_{ij},\qquad    i,j  =1,\dots ,n.
\end{gather*}  Moreover, the $n^2$ independent components $a_{ij}$
and $b_{ij}$ of $A$ in the basis~$\mathcal{B}$ can be obtained from
the components~$\tilde{a}_{ij}$ in the basis~$\widetilde{\mathcal{B}}$.

Now, for each $j\in\{1,\dots ,n\}$ we have that conditions
$i_{e_1\dots e_j}A=0$ give $jn$ independent restrictions on the
$n^2$-dimensional space $g^1_u(P_1)$. This implies that
\begin{gather*}
 \dim g^1_u(P_1)_{e_1\dots e_j}=n(n-j), \qquad  \dim g^1_u(P_1)_{e_1\dots e_n,v_1,\dots v_j}=0.
\end{gather*}
It follows that
\begin{gather*}
  \dim g^1_u\left(P_1\right)   + \sum_{i=1}^{n} \dim
  g^1_u\left(P_1\right)_{e_1,\dots ,e_i} + \sum_{i=1}^{n} \dim
  g^1_u\left(P_1\right)_{e_1,\dots ,e_n,v_1,\dots ,v_i}
  \\
 \qquad{}  = n^2 + n(n-1) + \cdots + n= n^2(n+1)/2=\dim g^2_u\left(P_1\right),
\end{gather*}
which shows that formula
\eqref{quasir} is satisf\/ied for $P_1$, $k=1$, and the basis
$\widetilde{\mathcal B}$. Therefore, $\widetilde{\mathcal B}$ is a~quasi-regular basis and hence the symbol
of~$P_1$ is involutive.
\end{proof}

\subsection{First obstruction for the projective
metrizability problem} \label{subsec:obstr1}

We have seen in the previous subsection that one condition, of the two
suf\/f\/icient conditions of the Cartan--K\"ahler theorem, for the formal
integrability of $P_1$, is satisf\/ied. In this subsection we address the
second suf\/f\/icient condition. We prove that there is only one obstruction for the formal
integrability of the projective metrizability operator $P_1$ and this is due to the
the curvature tensor $R$ of the induced nonlinear connection.

\begin{Theorem} \label{thm:obstr1} A first-order formal solution $\theta \in \Lambda^1_v$ of
the system \eqref{cjh} can be lifted into a second-order solution,
which means that $\overline{\pi}_{1}: R^{2}(P_1) \to R^{1}(P_1)$
is onto, if and only if
\begin{gather} d_R\theta=0, \label{drtheta} \end{gather}
where $R$ is the curvature tensor \eqref{curvature}.
\end{Theorem}
\begin{proof}
Using the notations from Subsection \ref{subsec:inv}, we denote by $K$, the cokernel of the morphism~$\sigma^2(P_1)$,
\begin{gather} K=\frac{T^*\otimes \left(T^*_v \oplus
\Lambda^2T^*_v \oplus \Lambda^2 T^*_v\right)}{\operatorname{Im}
\sigma^2(P_1)}. \label{spacek} \end{gather}
We will prove the theorem by using the following classical result of
homological algebra, see \cite[Proposition~1.1]{grifone00}. There
exists a morphism $\varphi: R^1(P_1) \to K$ such that the sequence
\[
 R^2(P_1)\overset{\overline{\pi}_1}{\longrightarrow} R^1(P_1)
\overset{\varphi}{\longrightarrow} K
\] is exact. In particular, the
morphism  $\overline{\pi}_1$ is onto if and only if~$\varphi=0$.

We will build the morphism $\varphi$ and show that for $\theta \in
\Lambda^1_v$ such that $j^1_u\theta \in R^1_u(P_1)\subset J^1_uT^*_v$, a~f\/irst-order
solution of~$P_1$ at $u\in TM\setminus\{0\}$ we have that
$\varphi_u\theta=0$ if and only if~$(d_R\theta)_u=0$. The morphism
$\varphi$ is represented in the diagram~\eqref{diagram} by dashed
arrows.

To build $\varphi$, we have to def\/ine f\/irst a morphism of vector bundles
\[
 \tau: \ T^*\otimes \left(T^*_v \oplus
\Lambda^2T^*_v \oplus \Lambda^2 T^*_v\right) \to K,
\] such that the
f\/irst row in the following diagram is exact.
\begin{gather}
  \label{diagram}
  {\diagram
 &  0 \dto  & 0 \dto & 0 \dto
  \\
 0\rto & g^2(P_1) \rto \dto & S^{2}T^*\otimes T^*_v \rto^{\sigma^2(P_1)} \dto^{\varepsilon} & T^* \otimes
  F_1 \rto \rdashed<1ex>|>{\tip}^{\tau} \dto^{\varepsilon} &
  K \rto & 0
  \\
 0 \rto &   R^{2}(P_1) \rto^{i} \dto^{\overline {\pi}_{1}} & J_{2} T^*_v \rto \dto^{\pi_{1}}
  \rdashed<1ex>|>{\tip}^{p^1(P_1)}
  & J^1 F_1 \dto^{\pi} \udashed<1ex>|>{\tip}^{\nabla}
  \\
  0 \rto & R^1(P_1) \rto\rdashed<1ex>|>{\tip}^i  & J^1T^*_v \rto^{p^o(P_1)}
  \udashed<1ex>|>{\tip}  \dto & F_1 \dto
  \\
&  & 0 & 0
  \enddiagram}
\end{gather}

For the vector bundle $K$ given by formula \eqref{spacek}, the
dimension of its f\/ibres is $n^2(n-1)/2$. Therefore, we can view this
vector bundle over $TM\setminus\{0\}$ as follows
\[
K=\oplus^{(2)}\Lambda^2T^*_v \oplus^{(3)}\Lambda^3T^*_v.
\]
Therefore the map $\tau$ has $5$ components.
The f\/ive components of $\tau=(\tau_1,\dots , \tau_5)$, are given as
follows
\begin{alignat*}{4}
& \tau_{1}(A,B_1,B_2) = \tau_JA - i_{\mathbb{C}}B_1, \qquad &&
\tau_{2}(A,B_1,B_2) = \tau_hA - i_{\mathbb{C}}B_2, \qquad &&
\tau_{3}(A,B_1,B_2) = \tau_JB_1,&  \\
& \tau_{4}(A,B_1,B_2) = \tau_hB_2, \qquad &&
\tau_{5}(A,B_1,B_2) = \tau_hB_1+\tau_JB_2, &&&
\end{alignat*}
for $A\in T^*\otimes T^*_v$, $B_1, B_2 \in T^*\otimes
\Lambda^2T^*_v$. Using the above def\/inition of the f\/ive components of $\tau$,
formula \eqref{sigma2p1} that def\/ines the three components of
$\sigma^2(P_1)$, and the symmetry in the f\/irst two arguments of an
element $B\in S^2T^*\otimes T^*_v$ we can prove $\left(\tau\circ
\sigma^2(P_1)\right)(B)=0$. For example, the f\/irst component of this
composition is given by {\sloppy
\begin{gather*}
\left(\tau_1\circ \sigma^2(P_1)\right)(B)(X,Y)  =
\left(\tau_J\sigma^2({\mathcal{L}_{\mathbb{C}}})\right)(B)(X,Y)-
\left(i_{\mathbb{C}}\sigma^2(d_J)\right)(B)(X,Y)\\
 \phantom{\left(\tau_1\circ \sigma^2(P_1)\right)(B)(X,Y)}{}
 = B(JX, \mathbb{C}, Y)- B(JY, \mathbb{C}, X) - B(\mathbb{C}, JX, Y)+ B(\mathbb{C}, JY, X)=0.\!
\end{gather*} It follows that $\operatorname{Im} \sigma^2(P_1)
\subset \operatorname{Ker} \tau $. By comparing the dimensions, it
is easy to see that  \mbox{$\operatorname{Im} \sigma^2(P_1) =
\operatorname{Ker} \tau $}, and therefore the f\/irst row in diagram~\eqref{diagram} is exact.

}

Consider $\nabla$ a linear connection on $TM\setminus\{0\}$ such
that~$\nabla J=0$. It follows that the connection~$\nabla$ preserves
the vertical distribution and hence it will preserve semi-basic
forms. Therefore, one can view $\nabla$ as a connection on the
f\/ibre bundle $F_1\to TM\setminus\{0\}$.  Using Lemma~\ref{lem:nablaomega}, it follows that derivations
$\mathcal{D}_{\mathbb{C}}=i_{\mathbb{C}}\nabla$,
$\mathcal{D}_{J}=\tau_{J}\nabla$, and $\mathcal{D}_{h}=\tau_{h}\nabla$
preserve semi-basic forms. As a~f\/irst-order partial
dif\/ferential operator, we can identify connection $\nabla$ with
the bundle morphism $p^0(\nabla): J^1F_1 \to T^*\otimes F_1$.
We will use this bundle morphism to def\/ine the map $\varphi: R^1(P_1)
\to K$ we mentioned at the beginning of the proof.

Consider $\theta\in \Lambda^1_v$ such that $j^1_u\theta\in
R^1_u(P_1)\subset J_u^1T^*_v$ is a f\/irst-order solution of $P_1$
at $u\in TM\setminus\{0\}$. Then, we def\/ine
\[
 \varphi_u\theta=\tau_u\nabla P_1\theta=\tau_u(\nabla \mathcal{L}_{\mathbb{C}}\theta,
\nabla d_J\theta, \nabla d_h\theta).
\] We will compute now the
f\/ive components of map $\varphi$. Since
$\mathcal{L}_{\mathbb{C}}\theta$, $d_J\theta$, and $d_h\theta$
vanish at $u\in TM\setminus\{0\}$, using Lemma \ref{lem:dec_DL},
it follows that when acting on this semi-basic forms we have
$\mathcal{D}_{\mathbb{C}}=\mathcal{L}_{\mathbb{C}}$,
$\mathcal{D}_{J}=d_J$, and $\mathcal{D}_{h}=d_h$. Using the fact that
$[J, \mathbb{C}]=J$, $[h, \mathbb{C}]=0$,  $[J, J]=0$, and $[h,J]=0$,
it follows that
\begin{gather*}
\tau_{1}\left(\nabla P_1\theta\right)_u  =  \left(\tau_J\nabla
\mathcal{L}_{\mathbb{C}} \theta - i_{\mathbb{C}}\nabla
d_J\theta\right)_u= \left(d_J\mathcal{L}_{\mathbb{C}} \theta -
\mathcal{L}_{\mathbb{C}} d_J\theta\right)_u =
(d_{[J,\mathbb{C}]}\theta)_u=0; \\
\tau_{2}\left(\nabla P_1\theta\right)_u  =  \left(\tau_h\nabla
\mathcal{L}_{\mathbb{C}} \theta - i_{\mathbb{C}}\nabla
d_h\theta\right)_u= \left(d_h\mathcal{L}_{\mathbb{C}} \theta -
\mathcal{L}_{\mathbb{C}} d_h\theta\right)_u =
(d_{[h,\mathbb{C}]}\theta)_u=0; \\
\tau_{3}\left(\nabla P_1\theta\right)_u  =  \left(\tau_J\nabla
d_J\theta\right)_u= \left(d^2_J\theta\right)_u
=\frac{1}{2}(d_{[J,J]}\theta)_u=0; \\
\tau_{4}\left(\nabla P_1\theta\right)_u  =  \left(\tau_h\nabla
d_h\theta\right)_u= \left(d^2_h\theta\right)_u
=\frac{1}{2}(d_{[h,h]}\theta)_u=(d_{R}\theta)_u; \\
\tau_{5}\left(\nabla P_1\theta\right)_u  =  \left(\tau_h\nabla
d_J\theta + \tau_J\nabla d_h\theta\right)_u=
(d_{[h,J]}\theta)_u=0.
\end{gather*}
From the above calculations it follows that a f\/irst-order formal
solution $\theta$ of the system~\eqref{cjh} can be lifted into a
second-order solution if and only if $d_R\theta=0$.
\end{proof}

Using notation \eqref{fdt} we can rewrite obstruction condition~\eqref{drtheta} as an algebraic Bianchi identity for the curvature
tensor
\begin{gather}
Fi_Rd\theta=h_{ik}R^k_{jl}+ h_{lk}R^k_{ij}+h_{jk}R^k_{li}=0. \label{Bianchi}
\end{gather}
An alternative expression for the algebraic Bianchi identity
\eqref{Bianchi} was obtained by Szilasi and Vattam\'any in
\cite[4.9.1a]{szilasi02}.

Using formula \eqref{curvature}, we obtain that any
solution of the system \eqref{cjh} necessarily satisf\/ies the curvature
obstruction~\eqref{drtheta}. In the next section we will discuss various cases when the obstruction~\eqref{drtheta} is automatically satisf\/ied. Another possibility, which we leave for
further work, is to add this obstruction to the projective metrizability
operator~$P_1$. In this case we can consider the f\/irst-order partial dif\/ferential operator~$P_2:
\Lambda^1_v(TM\setminus\{0\}) \to \Lambda^1_v(TM\setminus\{0\})
\oplus \Lambda^2_v(TM\setminus\{0\}) \oplus
\Lambda^2_v(TM\setminus\{0\}) \oplus
\Lambda^3_v(TM\setminus\{0\}) $,
\begin{gather*}
P_2  = \left({\mathcal L}_{\mathbb C}, d_J, d_h, d_R \right).
\end{gather*}
Following a similar approach as we did for the projective
metrizability operator~$P_1$, we can use the Cartan--K\"ahler theorem to study
the formal integrability of the dif\/ferential operator $P_2$.

\section{Classes of sprays that are projectively metrizable}
\label{sec:cpm}

In this section we present three classes of sprays for which the
projective metrizability opera\-tor~$P_1$ is formally integrable and hence
the system \eqref{cjh} always has solutions. Therefore, for each of
these classes we address the projective metrizability problem, by
discussing the set of algebraic conditions \eqref{palgebric} only,
which as we show are always satisf\/ied. We
will also provide examples of projectively metrizable sprays that are
not Finsler metrizable.

\subsection{Projectively metrizable sprays} \label{subsec:pms}
In this subsection we assume that a spray $S$ is analytic, on an
analytic manifold $M$. We show that if for spray $S$ the
projective metrizability operator $P_1$ is formally integrable then the
spray is projectively metrizable.

For a semi-basic $1$-form $\theta=\theta_i(x, y)dx^i \in \Lambda^1_v$,
we will express its f\/irst-order jet $j^1\theta\in J^1T^*_v$ in the
adapted dual basis $\{dx^i, \delta y^i\}$, induced by the nonlinear
connection associated to the spray, which means
\[
j^1\theta=\frac{\delta \theta_i}{\delta x^j} dx^j\otimes dx^i +
\frac{\partial \theta_i}{\partial y^j}\delta y^j\otimes dx^i.
\] This
expression provides us local coordinates $(x^i, y^i, \theta_i,
\theta_{ij}, \theta_{i\underline{j}})$ for $J^1T^*_v$. The typical
f\/ibre for the f\/ibre bundle $J^1T^*_v \to TM\setminus\{0\}$ is
$\mathbb{R}^{n*}\times L_2(n,\mathbb{R}) \times
L_2(n,\mathbb{R})$. With respect to these local coordinates, the f\/ibre
$R_u^1(P_1)$ of f\/irst-order formal solution of $P_1$ at $u=(x^i,y^i)\in
TM\setminus\{0\}$ can be expressed as follows
\[
 R^1_u(P_1)=\{(x^i, y^i, \theta_i,
\theta_{ij}, \theta_{i\underline{j}}) \in J^1_uT^*_v,
\theta_{ij}=\theta_{ji}, \theta_{i\underline{j}}=
\theta_{j\underline{i}}, \theta_{i\underline{j}}y^j=0\}.
\] Hence the
typical f\/ibre of the f\/ibre bundle $R^1(P_1) \to TM\setminus\{0\}$ is
$\mathbb{R}^{n*}\times L_{2, s}(n,\mathbb{R}) \times
L_{2,s}(n-1,\mathbb{R})$, where $L_{2,s}(n,\mathbb{R})$ is the space
of bilinear symmetric forms on $\mathbb{R}^n$.

Consider $\theta$ a solution of the system \eqref{cjh}, with the initial
data $(\theta^0_i, \theta^0_{ij}, \theta^0_{i\underline{j}})\in R^1_u(P_1)$ satisfying
the algebraic conditions~\eqref{palgebric}. This means $\theta^0_iy^i>0$
(in each f\/ibre, $y^i$ is a f\/ixed direction, hence one can choose
$(\theta^0_i)\in \mathbb{R}^{n*}$ such that $\theta^0_iy^i>0$) and
rank$(\theta^0_{i\underline{j}})=n-1$ (choose
$(\theta^0_{i\underline{j}})\in L_{2,s}(n-1, \mathbb{R})$ and extend
it to $\mathbb{R}^{n-1}\oplus \{a(y^i), a\in \mathbb{R}\}$ such that
$\theta^0_{i\underline{j}}y^j=0$).  If we assume that $M$ is
connected and $\dim M\geq 2$, then $TM\setminus\{0\}$ is also connected.
Therefore, due to continuity, the solution $\theta$ satisf\/ies the algebraic
conditions~\eqref{palgebric}, on the connected component of $u\in TM\setminus\{0\}$ where
$\theta$ is def\/ined.

We present now some classes of sprays for which the projective
metrizability operator $P_1$ is always integrable, and hence these
sprays will be projectively metrizable. These classes of sprays are:
\begin{itemize}\itemsep=0pt
\item[i)] f\/lat sprays, $R=0$;
\item[ii)] isotropic sprays, $R=\alpha\wedge J+\beta\otimes \mathbb{C}$, for $\alpha$ a
semi-basic $1$-form and $\beta$ a semi-basic $2$-form on $TM\setminus\{0\}$;
\item[iii)] arbitrary sprays on $2$-dimensional manifolds.
\end{itemize}
For each of these classes of sprays, we will show that the curvature obstruction
is automatically satisf\/ied and hence the projective metrizability
problem will always have a solution in the analytic case.

In the f\/lat case, the obstruction is automatically satisf\/ied. The fact
that f\/lat sprays are projectively metrizable was already demonstrated
with other methods in \cite{crampin07a}.

Assume that a spray $S$ is isotropic. It follows that the curvature tensor has the form  $R=\alpha\wedge
J+\beta\otimes \mathbb{C}$,  for $\alpha\in \Lambda^1_v$ and $\beta\in
\Lambda^2_v$.  Then, for a semi-basic $1$-form $\theta$ on
$TM\setminus\{0\}$, we have
\begin{gather} d_R\theta=\alpha\wedge d_J\theta + \beta \otimes
  \mathcal{L}_{\mathbb{C}}\theta. \label{eq:driso}\end{gather}
If $\theta$ is a solution of the dif\/ferential system~\eqref{cjh} it follows that~$\mathcal{L}_{\mathbb{C}}\theta=0$ and
$d_J\theta=0$, and using formula~\eqref{eq:driso} it follows that
$d_R\theta=0$. Therefore, the obstruction for the formal integrability
of~$P_1$ is satisf\/ied. In~\cite{crampin07b} it has been shown that any
isotropic sprays is projectively equivalent to a f\/lat spray and hence
it is projectively metrizable.

If $\dim M=2$ then  for a semi-basic $1$-form $\theta$ on
$TM\setminus\{0\}$,  $d_R\theta$ is a semi-basic $3$-form and hence it
will have to vanish. It has been shown by Matsumoto~\cite{matsumoto95}
that every spray on a surface is projectively related to a Finsler
space, using the original discussion of Darboux~\cite{darboux94} about
second-order dif\/ferential equations.

\subsection{Examples}
In this subsection we provide examples of non-metrizable Finsler
sprays in the last two of the above mentioned classes of projectively
metrizable sprays.

Consider the following system of second-order ordinary  dif\/ferential equations
in some open domain in $\mathbb{R}^2$, which was proposed by Anderson
and Thompson in \cite[Example~7.2]{anderso92}:
\begin{gather}
\frac{d^2x^1}{dt^2} + \left(\frac{dx^1}{dt}\right)^2+
\left(\frac{dx^2}{dt}\right)^2 = 0, \qquad
\frac{d^2x^2}{dt^2} + 4\frac{dx^1}{dt} \frac{dx^2}{dt} =
0. \label{ex:sodeat}\end{gather}
The corresponding spray is
\begin{gather}
S=y^1\frac{\partial}{\partial x^1} + y^2\frac{\partial}{\partial x^2}
-\left((y^1)^2+ (y^2)^2\right) \frac{\partial}{\partial y^1} - 4y^1y^2
\frac{\partial}{\partial y^2}. \label{ex:sat} \end{gather}
It has been shown in \cite{anderso92} that the system~\eqref{ex:sodeat} is not variational and therefore the corresponding
spray~$S$ in formula \eqref{ex:sat} is not Finsler metrizable.  We can also use the techniques
from~\cite{muzsnay06} to show that the spray $S$ in formula~\eqref{ex:sat} is not
Finsler metrizable. However, according to the discussion in the
previous subsection, the spray $S$ is projectively metrizable.

Next we consider another example of projectively metrizable spray that
is not Finsler metrizable, which was proposed by G.~Yang in~\cite{yang11}. Consider~$F$ a projectively f\/lat Finsler function on some open domain $U\subset
\mathbb{R}^n$ \cite[\S~13.5]{shen01}. This means that the geodesic
spray $S$ of $F$ is projectively equivalent to a f\/lat
spray. Therefore, spray $S$ is locally given by:
\begin{gather*}
S=y^i\frac{\partial}{\partial x^i}- 2P(x,y)y^i\frac{\partial}{\partial
  y^i}, 
  \end{gather*} where $P$ is $1$-homogeneous
function on $U\times (\mathbb{R}^n\setminus\{0\})$.
We assume that for the projectively f\/lat Finsler function $F,$ its f\/lag curvature is
constant $\kappa \in \mathbb{R}$, $\kappa \neq 0$, \cite[\S~3.5]{grifone00}, \cite[\S~11.1]{shen01}.
This is equivalent to the fact that the Jacobi endomorphism induced by
the spray $S$ has the form
\begin{gather*}
\Phi=\kappa\big(F^2J - Fd_JF\otimes \mathbb{C}\big). 
\end{gather*}
Yang shows in \cite[Theorem~1.2]{yang11} that the projective metrizability class of
$S$ contains sprays that are not Finsler metrizable. More precisely, he
shows that for $\lambda \in \mathbb{R}$ such that $\lambda\neq 0$ and
$\kappa+\lambda^2\neq 0$,  then the spray
\begin{gather} \widetilde{S}=S-2\lambda
  F\mathbb{C}\label{eq:splf} \end{gather} cannot be projectively
f\/lat and hence it is not Finsler metrizable.

For spray $\widetilde{S}$ one can compute the corresponding
geometric structures: nonlinear connection,
Jacobi endomorphism, curvature tensor in terms of the corresponding
ones induced by spray $S$:
\begin{gather}
\widetilde{h}   =   h+[P\mathbb{C}, J], \qquad
\label{eq:stildes} \widetilde{\Phi}   =   \Phi+ \lambda^2(F^2J- F d_JF \otimes
\mathbb{C}), \qquad
\widetilde{R}   =   R + \lambda^2 Fd_JF \wedge J.
\end{gather}
Therefore $\widetilde{S}$ has constant f\/lag curvature $\kappa+\lambda^2$
and it is also isotropic. Then one can also use formulae~\eqref{eq:stildes} and Theorem~2 from \cite{muzsnay06}, or Theorem~7.2 from~\cite{grifone00}, to show
that Yang's example  given in formula~\eqref{eq:splf}  is not Finsler
metrizable. Yang's example can be extended and it can be shown that
for an arbitrary spray, its projective class contains sprays that are not Finsler
metrizable,~\cite{bm11}.

Therefore, spray $\widetilde{S}$ in formula \eqref{eq:splf}  is
projectively metrizable but it is not Finsler metrizable.

\subsection*{Acknowledgements}
The work of IB was supported by the Romanian National
Authority for Scientif\/ic Research, CNCS UEFISCDI, project number
PN-II-RU-TE-2011-3-0017. The work of Z.M. has been
supported by the Hungarian Scientif\/ic Research Fund (OTKA) Grant K67617.

\pdfbookmark[1]{References}{ref}
\LastPageEnding


\begin{thebibliography}{99}

\footnotesize\itemsep=0pt

\bibitem{paiva05}
\'Alvarez Paiva J.C.,
Symplectic geometry and Hilbert's fourth problem,
{\it J.~Differential Geom.} {\bf 69} (2005), 353--378.

\bibitem{anderso92}
Anderson I., Thompson G.,
The inverse problem of the calculus of variations for ordinary dif\/ferential equations,
{\it Mem. Amer. Math. Soc.} {\bf 98} (1992),  no.~473.

\bibitem{aim93}
Antonelli P.L., Ingarden  R.S., Matsumoto M.,
The theory of sprays and Finsler spaces with applications in physics and biology,
Kluwer Academic Publisher, Dordrecht, 1993.

\bibitem{bcd10}
Bucataru I., Constantinescu O.A., Dahl M.F.,
A geometric setting for systems of ordinary dif\/ferential equations,
\href{http://dx.doi.org/10.1142/S0219887811005701}{{\it Int. J. Geom. Methods Mod. Phys.}} {\bf 8} (2011), 1291--1327,
\href{http://arxiv.org/abs/1011.5799}{arXiv:1011.5799}.

\bibitem{bucataru09}
Bucataru I., Dahl M.F.,
Semi basic 1-forms and Helmholtz conditions for the inverse problem of the calculus of variations,
\href{http://dx.doi.org/10.3934/jgm.2009.1.159}{{\it J.~Geom. Mech.}}  {\bf 1} (2009), 159--180,
\href{http://arxiv.org/abs/0903.1169}{arXiv:0903.1169}.

\bibitem{bm11}
Bucataru I., Muzsnay Z.,
Projective and Finsler metrizability: parameterization-rigidity of the geodesics,
\href{http://arxiv.org/abs/1108.4628}{arXiv:1108.4628}.

\bibitem{BCG91}
Bryant R.L., Chern S.S.,  Gardner R.B., Goldschmidt H.L., Grif\/f\/its P.A.,
Exterior dif\/ferential systems, {\it Mathematical Sciences Research Institute Publications}, Vol.~18, Springer-Verlag, New York, 1991.

\bibitem{crampin81}
Crampin M.,
On the dif\/ferential geometry  of the Euler--Lagrange equation and the inverse problem of Lagrangian dynamics,
\href{http://dx.doi.org/10.1088/0305-4470/14/10/012}{{\it  J.~Phys.~A: Math. Gen.}} {\bf 14} (1981), 2567--2575.

\bibitem{crampin07a}
Crampin M.,
On the inverse problem for sprays,
{\it Publ. Math. Debrecen} {\bf 70}  (2007), 319--335.

\bibitem{crampin07b}
Crampin M.,
Isotropic and $R$-f\/lat sprays,
{\it Houston~J. Math.}  {\bf 33} (2007), 451--459.

\bibitem{crampin08}
Crampin M.,
Some remarks on the Finslerian version of Hilbert's fourth problem,
{\it Houston~J. Math.}   {\bf 37} (2011), 369--391.

\bibitem{darboux94}
Darboux G.,
Le\c cons sur la theorie des surfaces, III,
Gauthier-Villars, Paris, 1894.

\bibitem{frolicher56}
Fr\"olicher A., Nijenhuis A.,
Theory ot vector-valued dif\/ferential forms. I.~Derivations in the graded ring of dif\/ferential forms,
{\it Nederl. Akad. Wet. Proc. Ser.~A} {\bf 59} (1956), 338--359.

\bibitem{grifone72}
Grifone J.,
Structure presque-tangente et connexions.~I,
{\it Ann. Inst. Fourier (Grenoble)} {\bf 22}  (1972), 287--334.

\bibitem{grifone00}
Grifone J., Muzsnay Z.,
Variational principles for second-order dif\/ferential equations. Application of the Spencer theory to characterize variational sprays,
\href{http://dx.doi.org/10.1142/9789812813596}{World Scientif\/ic Publishing Co., Inc.}, River Edge, NJ, 2000.

\bibitem{hamel03}
Hamel G.,
\"Uber die Geometrien, in denen die Geraden die K\"urzesten sind,
\href{http://dx.doi.org/10.1007/BF01444348}{{\it Math. Ann.}} {\bf 57} (1903), 231--264.

\bibitem{klein68}
Klein J.,  Voutier A.,
Formes ext\'erieures g\'en\'eratrices de sprays,
{\it Ann. Inst. Fourier (Grenoble)} {\bf 18} (1968), 241--260.

\bibitem{KMS93}
Kol\'ar I., Michor P.W., Slovak J.,
Natural operations in dif\/ferential geometry,
Springer-Verlag, Berlin, 1993.

\bibitem{krupkova07}
Krupkov\'a O., Prince G.E.,
Second order ordinary dif\/ferential equations in jet bundles and the inverse  problem of the calculus of variations,
 in \href{http://dx.doi.org/10.1016/B978-044452833-9.50017-6}{Handbook of Global Analysis}, Editors D.~Krupka and D.J. Saunders,   Elsevier Sci. B.V., Amsterdam,  2007, 837--904.

\bibitem{deleon89}
de Le\'on M., Rodrigues P.R.,
Methods of dif\/ferential geometry in analytical mechanics,
{\it North-Holland Mathematics Studies}, Vol.~158, North-Holland Publishing Co., Amsterdam, 1989.

\bibitem{lovas07}
Lovas R.L.,
A note on Finsler--Minkowski norms,
{\it Houston J. Math.}  {\bf 33} (2007), 701--707.

\bibitem{matsumoto95}
Matsumoto M.,
Every path space of dimension two is projectively related to a Finsler space,
\href{http://dx.doi.org/10.1007/BF02228993}{{\it Open Syst. Inf. Dyn.}} {\bf 3} (1995), 291--303.

\bibitem{matsumoto86}
Matsumoto M.,
Foundations of Finsler geometry and special Finsler spaces, Kaiseisha Press, Shigaken, 1986.

\bibitem{miron94}
Miron R., Anastasiei M.,
The geometry of Lagrange spaces: theory and applications,
{\it Fundamental Theories of Physics}, Vol.~59, Kluwer Academic Publishers Group, Dordrecht, 1994.

\bibitem{morandi90}
Morandi G., Ferrario C., Lo Vecchio G., Marmo G., Rubano C.,
The inverse problem in the calculus of variations and the geometry of the tangent bundle,
\href{http://dx.doi.org/10.1016/0370-1573(90)90137-Q}{{\it Phys. Rep.}} {\bf 188} (1990), 147--284.

\bibitem{muzsnay06}
Muzsnay Z.,
The Euler--Lagrange PDE and Finsler metrizability,
{\it Houston J. Math.}  \textbf{32} (2006), 79--98,
\href{http://arxiv.org/abs/math.DG/0602383}{math.DG/0602383}.

\bibitem{rapcsak62}
Rapcs\'ak A.,
Die Bestimmung der Grundfunktionen projektiv-ebener metrischer R\"aume,
{\it Publ. Math. Debrecen} {\bf 9} (1962), 164--167.

\bibitem{sarlet82}
Sarlet W.,
The Helmholtz conditions revisited. A new approach to the inverse problem of Lagrangian dynamics,
\href{http://dx.doi.org/10.1088/0305-4470/15/5/013}{{\it J.~Phys.~A: Math. Gen.}} \textbf{15} (1982), 1503--1517.

\bibitem{shen01}
Shen Z.,
Dif\/ferential geometry of spray and Finsler spaces,
Kluwer Academic Publishers, Dordrecht, 2001.

\bibitem{szilasi07}
Szilasi J.,
Calculus along the tangent bundle projection and projective metrizability,
in \href{http://dx.doi.org/10.1142/9789812790613_0045}{Dif\/ferential Geometry and its Applications},
 World Sci. Publ., Hackensack, NJ, 2008, 539--558.

\bibitem{szilasi02}
Szilasi J., Vattam\'any S.,
On the Finsler-metrizabilities of spray manifolds,
\href{http://dx.doi.org/10.1023/A:1014928103275}{{\it Period. Math. Hungar.}} {\bf 44} (2002), 81--100.

\bibitem{yang11}
Yang G.,
Some classes of sprays in projective spray geometry,
\href{http://dx.doi.org/10.1016/j.difgeo.2011.04.041}{{\it Differential Geom. Appl.}} {\bf 29} (2011), 606--614.

\end{thebibliography}
\end{document}